\documentclass[11pt,a4paper,twoside]{amsart}
\RequirePackage[l2tabu, orthodox]{nag}

\usepackage[utf8]{inputenc}
\usepackage[english]{babel}
\usepackage[T1]{fontenc}
\usepackage{microtype, fancyhdr, lmodern, xcolor}

\usepackage[a4paper, hmarginratio=1:1]{geometry}

\usepackage{amsfonts, amsmath,  amssymb, mathrsfs, amscd}
\numberwithin{equation}{section}    
\usepackage{faktor}
\allowdisplaybreaks

\usepackage{url, enumitem}
\usepackage[noadjust]{cite}

\usepackage{multirow,bigdelim, caption}
\usepackage{booktabs}

\usepackage{graphicx}
\graphicspath{{Pictures/}}

\usepackage{hyperref}

\newcommand{\todopm}[1]{{\textcolor{blue}{#1}}}
\newcommand{\ppm}[1]{{\textcolor{red}{#1}}}

\theoremstyle{plain}
\newtheorem{theorem}{Theorem}[section]
\newtheorem{proposition}[theorem]{Proposition}
\newtheorem{corollary}[theorem]{Corollary}
\newtheorem{lemma}[theorem]{Lemma}
\newtheorem{conjecture}[theorem]{Conjecture}

\theoremstyle{remark}
\newtheorem{remark}[theorem]{Remark}

\newtheorem{question}[theorem]{Question}

\theoremstyle{definition}
\newtheorem{definition}[theorem]{Definition}

\newcommand\Nn{\mathbb{N}}
\newcommand\Zz{\mathbb{Z}}
\newcommand\Qq{\mathbb{Q}}
\newcommand\Rr{\mathbb{R}}
\newcommand\Cc{\mathbb{C}}
\newcommand{\C}{\mathbb{C}}

\newcommand\BB[1]{B_{#1}}	
\newcommand{\B}{B}              
\newcommand{\rhob}{{\mathbf \varrho}}  

\newcommand\ii{i}
\newcommand\inv{^{-1}} 
\newcommand\jj{j}

\newcommand\LB[1]{LB_{#1}}	
\newcommand\LH{L\!H}	
\newcommand{\LBcat}{{\mathsf{ LB}}}  
\newcommand{\Bcat}{{\mathsf{ B}}}    
\newcommand{\Hcat}{{\mathsf{ H}}}    
\newcommand{\LHcat}{{{LH}}}  
\newcommand{\Mat}{{\mathsf{ Mat}}}   
\newcommand{\Vect}{{\mathsf{ Vect}}} 
\newcommand{\ann}{Ann\,}
\newcommand{\fb}[1]{\framebox{#1}} 
\newcommand{\SPe}{S\!P}
\newcommand{\Gx}{G_{\chi}} 
\newcommand{\func}[2]{#1-\!\!\!\!\mod \rightarrow #2\!-\!\!\!\!\!\mod} %
\newcommand{\SP}{Burau--Rittenberg}
\newcommand{\FF}{{\mathsf F}} 

\newcommand\nn{n}
\newcommand\nno{{n-1}}

\newcommand\rr[1]{\rho_{\hspace{-0.3ex}#1}^{\null}}	
\newcommand{\pp}{p}
\newcommand\rrp[1]{\pp_{\hspace{-0.3ex}#1}^{}}	
\newcommand\rrb[1]{\rhob_{\hspace{-0.3ex}#1}^{\null}}	

\newcommand\sig[1]{\sigma_{\hspace{-0.3ex}#1}^{\null}} 
\newcommand\sigg[2]{\sigma_{\hspace{-0.3ex}#1}^{#2}}	

\newcommand{\R}{{\mathfrak R}}  
\newcommand{\Q}{{\mathfrak Q}}  
\newcommand{\A}{{\mathfrak A}}  

\newcommand{\Scat}{{\mathsf S}} 
\newcommand{\MMM}{{\mathcal M}}  

\newcommand{\Mot}{{\mathcal Mo}}  
\newcommand{\un}[1]{\underline{#1}} 

\newcommand{\e}{{\mathbf e}}
 
\newcommand{\Ccat}{{\mathsf C}}

\newcommand{\beq}{\begin{equation}}
\newcommand{\eq}{\end{equation}} 
\newcommand{\mat}{\left( \begin{array}}
\newcommand{\tam}{\end{array} \right)} 

\newcommand{\ignore}[1]{}  

\newcommand{\ppmm}[1]{#1}  

\title{Generalisations of Hecke algebras from \\ Loop Braid Groups}
\author{Celeste Damiani,
  Paul Martin,
  Eric C. Rowell
}
\date{\today}

\begin{document}

\begin{abstract}
We introduce a generalisation $\LH_\nn$ of the ordinary Hecke algebras informed
by the loop braid group $LB_n$ and the extension of the Burau representation
thereto.
The ordinary Hecke algebra has many remarkable arithmetic and
representation theoretic properties, and many applications.
We show that $\LH_\nn$ has analogues of several of these properties.
In particular we 
\ppmm{consider}
a class of local (tensor space/functor)
representations of the braid group
derived from a meld of the (non-functor) Burau representation \cite{Burau:1935}
and the (functor) \ppmm{Deguchi {\em et al}-Kauffman--Saleur-}Rittenberg representations
\cite{deguchi_89,deguchi1989,Deguchi_1990,kauffman1991,Rittenberg92},
here 
called \SP\ representations.
In its most supersymmetric case
somewhat mystical cancellations of anomalies occur so that 
the \SP\ 
representation extends to a loop
\SP\ representation.
And this factors through $\LH_\nn$.
Let $\SPe_\nn$ denote the corresponding (not necessarily proper)
quotient algebra,
$k$ the ground ring,
and $t \in k$ the loop-Hecke parameter.
We prove the following:
\begin{enumerate}
\item \label{Claim:dimension} $\LH_\nn$ is finite dimensional
  over a field.
\item The natural inclusion $LB_n \hookrightarrow LB_{n+1}$  passes to
  an inclusion 
  $\SPe_\nn  \hookrightarrow \SPe_{n+1}$. 
\item \label{Claim:radical}
Over $k=\Cc$,  
$\SPe_\nn / rad $ is  generically
the sum of
  simple matrix algebras of dimension
  (and Bratteli diagram)
  given by Pascal's triangle.
  (Specifically $\SPe_\nn / rad \cong \Cc S_n / \e^1_{(2,2)}$ where
  $S_n$ is the symmetric group and $\e^1_{(2,2)}$ is a $\lambda=(2,2)$
  primitive idempotent.)
\item We 
  determine the other fundamental invariants of $\SPe_\nn$
  representation theory: the Cartan
  decomposition matrix; and the quiver, which is of type-A.
\item \label{Claim:dependance} the structure of
  $\SPe_\nn  $   
  is independent of the parameter $t$, except for
  $t= 1$.
\item For $t^2 \neq 1$ then $\LH_n \cong \SPe_n$ at least up to rank
  $n=7$ (for $t=-1$ they are not isomorphic for $n>2$;
  for $t=1$ they are not isomorphic for $n>1$).
\end{enumerate}
Finally we discuss a number of other intriguing points arising from this
construction in
topology, representation theory and combinatorics.
\end{abstract}

\maketitle


\section{Introduction}

Until the 1980s, 
methods to construct linear representations of the braid group
$\B_n$ were relatively scarce.  
We have those factoring through the symmetric group and
the Burau representation \cite{Burau:1935},
and those factoring through the Hecke algebra \cite{Hoefsmit74} and
the Temperley--Lieb algebra \cite{TemperleyLieb71};
and, 
as for every group,
the closure in the monoidal category $Rep(\B_n)$.
These proceed essentially through
`combinatorial' devices such as 
Artin's presentation.
Then there are some more intrinsically `topological' constructions 
such as Artin's representation \cite{Artin}
(and Burau can be recast in this light \cite{LongPaton:1993}).

In the 80s there were 
notable steps forward.
Algebraic formulations of the Yang--Baxter equation began to
yield representations (see e.g. \cite{Baxter82}). 
And
Jones' discovery \cite{Jones:Hecke} of
link invariants
from finite dimensional quotients of the group algebra
$\mathbb{K}[\B_n]$ 
inspired a
revolution in braid group representations and topological invariants
\cite{Kauffman:1990,Birman-Wenzl:1989,Murakami:1987,HOMFLY,Wenzl:1988}.
Work of Drinfeld, Reshetikhin, Turaev, Jimbo,
\cite{Drinfeld:1987,Reshetikhin-Turaev:1991,Jimbo:1986} and others on quantum groups
yielded yet further representations.
Enriched through modern category theory
\cite{Turaev:2010,Kassel-Turaev:libro,Bakalov-Kirillov:2001,Damiani2019},
constructions are now relatively abundant.

The connections among $\B_n$ representations,
$2+1$-dimension topological quantum field theory (see e.g. \cite{Witten:1989})
and
statistical mechanics
(see e.g. \cite{Baxter82,akutsu87,Martin88,deguchi_89,Deguchi_1990})  
were already well established in the 1980s.
Even more recently, the importance of such representations in topological phases of matter \cite{FreedmanKitaevLarsenWang,RowellWang} in two spacial dimensions has led to an invigoration of interest, typically focused on unitary representations associated with the $2$-dimensional part of a $(2+1)$-TQFT.  In this context the braid group is envisioned as the group of motions of point-like quasi-particles in a disk, with the trajectories of these \emph{anyons} forming the braids in $3$-dimensions.  Here the braid group generators $\sigma_i$ correspond to exchanging the positions of the $i$ and $i+1$st anyons. The density of such braid group representations in the group of (special) unitary matrices is intimately related to the universality of quantum computational models built on these topological phases of matter \cite{FKW,FLW}, as well as the (classical) computational complexity of the associated link invariants \cite{Rowellparadigms}.  Indeed, there is a circle of conjectures relating \emph{finite} braid group images \cite{Naidu-Rowell:2011,Rowell-Stong-Wang:2009}, \emph{classical} link invariants, \emph{non-universal} topological quantum computers and \emph{localizable} unitary braid group representations \cite{Rowell-Wang:2012,GHRIMRN}.  The other side of this conjectured coin relates the holy grail of universal topological quantum computation with powerful $3$-manifold invariants through surgery on links in the three sphere.

\medskip

What is a non-trivial generalization of the braid group to
$3$-dimensions? Natural candidates are groups of motions:
heuristically, the
elements
are classes of trajectories of a compact
submanifold $N$ inside an ambient manifold $M$ for which the initial
and final positions of $N$ are set-wise the same.
The group of
motions
of points in a $3$-manifold
in effect
simply permutes the points, but
the motion of circles or more general links in a $3$-manifold is
highly non-trivial.  This motivates the study of these $3$-dimensional
motion groups, as defined in the mid-20th century by Dahm
\cite{Dahm:1962} and expounded upon by Goldsmith
\cite{Goldsmith:MotionGroups,Goldsmith:1982}.

More formally, a motion of $N$ inside $M$ is an ambient isotopy
$f_t(x)$ of $N$ in $M$ so that $f_0=id_M$ and $f_1(N)=N$.
Such a motion is \emph{stationary} if $f_t(N)=N$ for all $t$;
and given any motion $f$, we have the usual notion of the reverse
$\bar{f}$.
We say two motions $f,g$ are equivalent if the composition of $f$ with
$\bar{g}$ (via concatenation) gives a motion
endpoint-fixed
homotopic to a stationary motion as isotopies
$M\times [0,1]\rightarrow M$.
The
\emph{motion group} $\Mot(M,N)$ is the group of
motions modulo 
this equivalence.
When
$M$ and $N$ are both oriented we will consider only motions $f$ so
that $f_1(N)=N$ as an oriented submanifold, although one may consider
the larger groups allowing for orientation reversing motions.   

The motion groups of links inside $\mathbb{R}^3, S^3$ or $D^3$ and
their representations
are 
very rich, and only recently explored in the literature \cite{Bellingeri-Bodin:2016,DamianiKamada,Kadar-Martin-Rowell-Wang:2016,Bullivant-et-al:2020,Baez-Wise-Crans:2007,BMM19}.  Further enticement is provided by the prospect of applications to $3$-dimensional topological phases of matter with loop-like excitations (i.e. vortices) \cite{LevinWang}.  The fruitful symbiosis between braid group representations and $2$-dimensional condensed matter systems give us hope that $3$-dimensional systems could enjoy a similar relationship with  motion group representations, $3+1$-TQFTs, and invariants of surfaces embedded in $4$-manifolds (see e.g. \cite{Kamada:Markov,CKS}).

There are a few hints in the literature that the $3+1$-dimensional story has some key differences from the $2+1$-dimensional situation.  Reutter \cite{Reutter:preprint} has shown that semisimple $3+1$-TQFTs cannot detect smooth structures on $4$-manifolds.  Wang and Qiu \cite{WangQiu:preprint} provided  evidence that the mapping class group and motion group representations associated with $3+1$-dimensional Dijkgraaf-Witten TQFTs are determined via dimension reduction by the corresponding $2+1$-dimensional DW theory.  As the representation theory of motion groups has been largely neglected until very recently, it is hard to speculate on precise statements analogous to the $2$-dimensional conjectures and theorems.  


\medskip

In this article we take hints from the
the classical works \cite{Burau:1935,Hoefsmit74},
from the braid group revolution \cite{Jones:1987},
and more directly from statistical mechanics
\cite{Deguchi_1990,kauffman1991,Rittenberg92,Deguchi92},
to study representations of the motion
group of free unlinked circles in $3$-dimensional space, the
\emph{loop braid group} $LB_n$.
Presentations of $LB_n$ are known 
(see \cite{Fenn-Rimanyi-Rourke:1997} and \cite{Damiani:Journey} and
references therein).
As $LB_n$ contains the braid group $\B_n$ as an abstract subgroup, a
natural approach to finding linear representations is to extend known
$B_n$ representations to $LB_n$.
This has been considered by various authors, see
e.g. \cite{Bruillard-Chang-Hong-Plavnik-al:2015,Bardakov:ExtendingRepresentations,Kadar-Martin-Rowell-Wang:2016}.
Another  idea is to look for finite dimensional quotients
of the group algebra,
mimicking the techniques of \cite{Jones:1987,Birman-Wenzl:1989}.
As non-trivial finite-dimensional quotients of the braid group are not
so easy to find,
we take a hybrid approach: we combine the extension of the Burau
representation to $LB_n$
\cite{Burau:1935,Bardakov:ExtendingRepresentations} with the Hecke
algebras $\mathcal{H}_n$ obtained from $\mathbb{Q}(t)[B_n]$ as the
quotient by the ideal generated by $(\sigma_i+1)(\sigma_i-t)$.
While the naive quotient of $\mathbb{Q}(t)[LB_n]$ by this ideal does not
provide a finite dimensional algebra, certain additional quadratic relations
(satisfied by the extended Burau representation) are sufficient for
finite dimensionality, with quotient denoted $\LH_n$.
We find a local representation of $\LH_n$ that aids in the analysis of
its structure --- the Loop Burau-Rittenberg representation.
One important feature of the algebras  $\LH_n$ is that they are not semisimple; in fact, the image of the loop Burau-Rittenberg representation has a $1$-dimensional center, but is far from simple.  Its semisimple quotient by the Jacobson radical gives an interesting tower of algebras with Bratteli diagram exactly Pascal's triangle. 


\medskip

Our results suggest new lines of investigation into motion group
representations.  What other finite dimensional quotients of motion
group algebras can we find (confer e.g. \cite{Banjo2013})?
What is the role of (non)-semisimplicity in such quotients?
Can useful topological invariants be derived from these quotients?
What do these results say about $3+1$-dimensional TQFTs?


%
{\em Outline of the paper:}
In Sec.\ref{S:Motivation} we recall the Burau representation
and corresponding  knot invariants.
In Sec.\ref{ss:loop} we introduce loop Hecke algebras and prove they are finite dimensional.
In Sec.\ref{ss:basic} we develop arithmetic tools (calculus) that we will need.
In Sec.\ref{ss:locreps} we construct our local representations
and hence prove our main structure Theorems.
In Sec.\ref{ss:lowrank} we apply the results from \S\ref{ss:locreps}
to $\LH_n$,
and
make several conjectures on the open cases with
$t^2 =1$.
We conclude with a discussion of new directions opened up
by this work. 

\medskip

\medskip

\noindent
{\bf Acknowledgements} 
\\
CD and PM would like to thank the Leverhulme Trust for support under
research project 
grant
RPG-2018-029:
``Emergent Physics From Lattice Models of Higher Gauge Theory'';
and
Jo\~{a}o Faria Martins (PI) for encouragement.
CD would like to thank Paolo Bellingeri and Emmanuel Wagner for their guidance during the time this idea was being born.
PM would like to thank 
Raphael Bennett-Tennenhaus 
and Paula Martin
for useful conversations.
We all thank Alex Bullivant (and hence EPSRC), 
and the University of Leeds,
for partial support of 
ER's visit to Leeds. ER is grateful for insightful comments from Vaughan Jones. ER is partially supported by US NSF grant DMS-1664359, a Presidential Impact Fellowship from Texas A\&M University and a Simons Foundation Fellowship.  This research was carried out while ER was participating in a semester-long program on Quantum Symmetries at MSRI, a research institute partially supported by US NSF grant DMS-1440140.

\section{Burau representation, Hecke algebra and 
  invariants of knots}
\label{S:Motivation}

Let $\un{n} := \; \{ 1,2,...,n \}$. Then the braid group $B_n$ may be
identified with
the motion group
$\Mot(\Rr^2 , \un{n} \times \{ 0 \})$.
Artin showed that, 
for $\nn \geq 1$,  
$\BB\nn$ admits the presentation 
\begin{equation} \label{eq:Br}
    \Bigg\langle
    \sig1, \ldots, \sig{n-1} 
    \Bigg\vert 
    \begin{aligned}
    \sig{i} \sig j &= \sig j \sig{i}  \, &\text{for } \vert  i-j\vert > 1 \\
    \sig{i} \sig {i+1} \sig{i} &= \sig{i+1} \sig{i} \sig{i+1} \, &\text{for } i=1, \dots, \nn-2 \\
    \end{aligned}
    \Bigg\rangle
\end{equation}
We will write $\A_n(\sig{})$ for the set of relations here.

We will also need the symmetric group $S_n$. In a `motion group
spirit' this can be identified with
$\Mot(\Rr^3 , \un{n} \times \{ 0 \} \times \{ 0 \} )$.
It can be presented as a quotient of $B_n$ by the relation ${\sig{1}}^2=1$
(however since we will often want to have both groups together we will
soon rename the $S_n$ generators).

\subsection{Burau representation}
\label{SS:Braid_Burau}

We define Burau representation $\rhob \colon \BB\nn \to GL_\nn(\Zz[t, t\inv])$ as
 follows:
\begin{equation}
\label{E:braid_burau}
    \sig\ii \mapsto I_{\ii-1} \oplus \begin{pmatrix} 1-t & t \\ 1 & 0 \end{pmatrix} \oplus I_{\nn-i-1}.
\end{equation}

The Burau representation
has  Jordan--Holder decomposition
into a $1$-dimensional representation
(the vector~$(1, \ldots, 1)^T$ remains fixed) and 
an $(\nn-1)$-dimensional irreducible representation 
known as 
reduced Burau representation $\bar{\rhob} \colon \BB\nn \to GL_\nno(\Zz[t, t\inv])$. 
\ignore{{
defined by:
\begin{equation}
    \sig\ii \mapsto I_{\ii-2} \oplus 
    \begin{pmatrix} 1 & -t & 0 \\ 0 & -t & 0 \\ 0 & -1 & 1 \end{pmatrix}
    \oplus I_{\nn-i-2}.
\end{equation}
\ppm{-what happens when $i=1$?}
}}
The decomposition is not split over $\Zz[t, t\inv]$ --- an inverse of
$t+1$ is needed (see later).

\begin{remark}
One can also use the transpose matrix of~\eqref{E:braid_burau}
(depending on orientation choices while building the ``carpark cover''
of the punctured disc in the homological definition of Burau).
The transpose fixes $(1, \ldots,1, t, t^2, 1, \ldots, 1)^T$. 
\end{remark}

\subsection{Facts about the Burau representation}
\label{SS:Burau_trivia}

\begin{enumerate}
\item Burau is unfaithful for $\nn \geq 5$ 
(Moody \cite{Moody:1991} proved unfaithfulness for $\nn \geq 9$, 
Long and Paton \cite{LongPaton:1993} for $\nn \geq 6$, Bigelow \cite{Bigelow:1999} for $\nn = 5$).
 
\item The case $\nn=4$ is open, Beridze and Traczyk \cite{BeridzeTraczyk:2018} recently published some 
advances toward closing the problem.
 
\item It is faithful for~$n=2,3$ (Magnus and Peluso \cite{MagnusPeluso:1969}).

\item If we consider the braid group in its mapping class group formulation, 
  it has a homological meaning 
  (attached \textit{a posteriori} 
  to it, since Burau 
  used only combinatorial aspects of matrices
  \cite{Burau:1935}).
  The Burau representation
  describes the action of braids on the first homology group of the (covering of) the punctured disk. On the other hand the Alexander polynomial is extracted from the presentation matrix of the first homology group of the knot complement (the Alexander matrix). When we close up a braid, each element of homology of the punctured disk on the bottom becomes identified with its image in the punctured disk at the top. At this point the Alexander matrix of the closed braid is (roughly) the Burau matrix of the braid with the modification of identifying the endpoints.

More specifically, let $K$ be a knot, and $b$ a braid such that 
$\hat{b}$ is equivalent to~$K$. 
Then the Alexander polynomial $\Delta_K(t)$ can be obtained by computing:
\[
    \Delta_K(t) = \frac{\det(\bar{\rhob}(b) - I_{n-1})}{1 + t + \ldots + t^\nno}.
\]

So one can think of the Alexander polynomial of $K \sim \hat{b}$ as a rescaling
of the characteristic polynomial of the image of $b$ 
in the reduced representation. 

\end{enumerate}


Representations of $B_n$ are partially characterised by the eigenvalue
spectrum of the
image of $\sig{i}$. 
Observe that  
\begin{equation}
\label{eq:loc}
  \rhob(  {\sig\ii}^2 ) = (1-t) \rhob({\sig\ii} ) + tI_\nn ,
\end{equation}
i.e., the eigenvalue spectrum is $Spec(\rhob(\sig{i}))= \{ 1, -t \}$. 
Recall also that Kronecker products obey
$Spec(A \otimes B) = Spec(A).Spec(B)$,
so
$Spec(\rhob(\sig{i})\otimes \rhob(\sig{i})
= \rhob\otimes\rhob(\sig{i})) = \{ 1, -t, t^2 \}$.
From this we see that the spectrum is fixed under tensor product only
if $t=\pm 1$
(cf. for example \cite{kauffman1991}). 
\ignore{{
\ppm{ [TO DO:  Perhaps explain a little how to compare and contrast with KS.
Left implicit here is that one can apply Schur functors to select
sub-representations
from the tensor powers --- the nominally easiest are the symmetric and
exterior powers. Since these are subreps they may select from the
enlarged
spectrum in a way that gets around the no-go theorem. This is
basically irrelevant to us, but addresses a point that a couple of
people raise when they see this work...
Alternatively it might be better just to add some notes of caution
(and another Deguchi ref) at
the appropriate point later - when we mention supersymmetry.]  }
}}

\newcommand{\sigT}{T_}  

\subsection{Hecke algebras}

Let $R$ be an integral domain and $q_1, q_2$ elements of~$R$ with $q_2$ invertible. We define the Hecke algebra $H^R_\nn(q_1, q_2)$ to be 
the algebra with generators
$\{1, T_1, \ldots, T_{\nn -1}\}$ 
and the following defining relations:
\begin{align}
  T_\ii T_\jj &= T_\jj T_\ii  \, 
  &\text{for } \vert  \ii-\jj\vert > 1 \\
  T_\ii T_{\ii+1} T_\ii &= T_{\ii+1} T_\ii T_{\ii+1} \, 
  &\text{for } \ii=1, \dots, \nn-2 \\
\label{e:quadratic_sigma}  T_\ii^2 &= (q_1 + q_2)T_\ii -q_1q_2 \,
\, &\text{for } \ii=1, \dots, \nn-1 .
\end{align}

\begin{remark} \label{rem22}
\begin{enumerate}

\item  Relation \eqref{e:quadratic_sigma} coincides with the 
characteristic equation of the images of the generators under the Burau
representation~\eqref{SS:Braid_Burau} when $(q_1, q_2) = (1, -t)$.
We denote the resulting $1$-parameter Iwahori-Hecke algebra by~$H^R_\nn(t)$.

\item If  $t=1$ then $H^R_\nn(t)$ is the group algebra~$R[S_\nn]$
  (the free $R$-module $R S_n$ made an $R$-algebra in the usual way).

\item \label{i:mapBn} There is a 
  map from $\BB\nn$ to $H^R_\nn(t)$ 
sending $\sig\ii$ to~$T_\ii$.
Thus representations of $H^R_\nn(t)$ are equivalent to 
representations of $\BB\nn$ for which the  generators satisfy
relation~\eqref{e:quadratic_sigma}.
This is described in~\cite[Section~3]{Bigelow:2006}, \cite[Section~4]{Jones:1987}
in \cite[\S5.7]{Martin:M1} and many other places.

\item Fixing~$R=\Cc$, point 
  \eqref{i:mapBn}
  allows us to think of $H^R_\nn(t)$ 
as being isomorphic to the quotient~$H_n(t) := \faktor{\Cc[\BB\nn]}{\sigg{\ii}{2} = (1-t)\sig\ii + t}$.

\item Using the map in~\eqref{i:mapBn} we can represent 
any element of $H_n(t)$ as a linear combination of braid diagrams. 
The quadratic relation can be seen as a \emph{skein relation} 
on elementary crossings. Knowing a basis for~$H_n(t)$ makes this fact usable.
\end{enumerate}
\end{remark}

\begin{question} Why these parameters and this quadratic relation?

  As noted,
  Hecke algebras can be defined with two units of $R$
  as parameters. 
  We chose to fix these parameters to $(1, -t)$
  because from this quotient one should  
  recover the Alexander
polynomial.
Choosing $(-1, t)$ one should get the quotient on which Ocneanu traces are defined (see \cite[Chapter~4.2]{Kassel-Turaev:libro}). 
With the Ocneanu trace being a $1$-parameter family over a $1$-parameter 
algebra, we end up with polynomials in two variables. 
These polynomials are attached to the 
braid diagrams
that we can see representing elements of $H_n(t)$. 
Moreover they are defined in such a way to respect Markov moves, 
so they are invariants for the closures of said braids. 
Hence, they are knot invariants.
The quadratic relation
from \ref{rem22}\eqref{i:mapBn} 
translates
the trace in a \emph{skein relation}. 
Through the Ocneanu trace (normalised) 
the invariant that is obtained is the HOMFLY-PT polynomial, 
which specialises in both Alexander and Jones. 
Each specialisation corresponds to factoring through a further quotient 
of the Hecke algebra
(in the case of Jones, this is
a quotient of 
the Temperley-Lieb algebra).
Below we
``reverse engineer'' this process. 
\end{question}

\newcommand{\I}{{\mathfrak I}}

\section{Generalising Burau  
  and Hecke to loop braid groups}
\label{ss:loop}
\subsection{The loop braid group}

Here $S^1$ denotes the unit circle.
We now consider the loop braid group
$$
\LB\nn =  \Mot(\Rr^3 , \un{n} \times S^1 )
$$
(see e.g.
\cite{Goldsmith:MotionGroups,Savushkina:1996,Fenn-Rimanyi-Rourke:1997,%
BrendleHatcher:2013,Damiani:Journey,Kadar-Martin-Rowell-Wang:2016,%
Bruillard-Chang-Hong-Plavnik-al:2015}).


\newcommand{\LBxn}{\langle \Xi_n | \Q_n \rangle}

Consider the
set $\Xi_n = \{  \sig\ii, \rr\ii , i=1,2,...,n-1 \}$ 
and group $\LBxn$
presented by generators $\sig\ii$ and $\rr\ii$,
and relations
$\Q_n$
as follows. 
The $\sig\ii$s obey the braid relations as in \eqref{eq:Br};
the $\rr\ii$s obey the braid relations and also
\beq \label{eq:pp1}
\rr\ii \rr\ii = 1
\eq
and then there are mixed braid relations
\beq \label{eq:Qm0}
\rr{i} \rr{i+1} \sig{i} = \sig{i+1} \rr{i} \rr{i+1},
\qquad
\eq
\beq \label{eq:Qm}
\rr{i} \sig{i+1} \sig{i} = \sig{i+1} \sig{i} \rr{i+1},
\eq
\beq \qquad\hspace{2.1cm}
\sig{i}\rr{i\pm j}=\rr{i\pm j}\sig{i} \qquad (j>1)
\qquad\qquad  (all \; distant \; commutators).
\eq

\newcommand{\cc}{1}
\newcommand{\bb}{{2}}

\begin{remark}
The first mixed relation \eqref{eq:Qm0} implies its reversed order counterpart: 
\beq \label{eq:Qm0r}
\sig{i} \rr{i+1} \rr{i} = \rr{i+1} \rr{i} \sig{i+1}
\eq
whereas the reversed order second mixed relation does not hold.
\\
The relations also imply
\beq \label{eq:M1}
{\rr\bb\sig\cc\rr\bb} =
    {\rr1\sig2\rr1}
\eq
\end{remark}


\bigskip

We have (see e.g. \cite{Fenn-Rimanyi-Rourke:1997}) that
\beq \label{th:present}
LB_n \;\cong\;  \LBxn .
\eq


It will be convenient to give an {\em algebra presentation} for the
group algebra.
Recall that 
in an algebra presentation inverses are not present automatically by
freeness, so we may put them in by hand as formal symbols and then
impose the inverse relations.
Thus as a presented algebra we have
\[
k \LBxn   \;=\;   \langle \Xi_n \cup \Xi_n^{-} \; | \; \Q_n, \I_n \rangle_k
\]
--- here $kG$ means the group $k$-algebra of group $G$; 
and $\langle -|- \rangle_k$ means a $k$-algebra presentation;
and $\I_n$ is the set of inverse relations $\sig{i}\sigma^-_i =1$.


\subsection{
The loop--Hecke algebra \texorpdfstring{$LH_n$}{LHn}}

With \S\ref{S:Motivation} in mind, there {\em is} a suitable generalisation of
the Burau representation to $LB_n$. 

\begin{proposition}[\cite{Vershinin:1996}] \label{pr:GBx}
The map on generators of $LB_n$ given by
\begin{equation}
    \sig\ii \mapsto I_{\ii-1} \oplus \begin{pmatrix} 1-t & t \\ 1 & 0 \end{pmatrix} \oplus I_{\nn-i-1}.
\end{equation}
\begin{equation}
    \rr\ii \mapsto I_{\ii-1} \oplus \begin{pmatrix} 0 & 1 \\ 1 & 0 \end{pmatrix} \oplus I_{\nn-i-1}.
\end{equation}
extends to a representation 
$\rrb{GB}: LB_n \rightarrow GL_n(\Zz[t,t^{-1}])$. 
\end{proposition}
\proof
Direct calculation.
\qed

This group representation is not faithful 
for~$\nn \geq 3$~\cite{Bardakov:ExtendingRepresentations},
and corresponds to an Alexander polynomial for welded knots.

\medskip \bigskip



\newcommand{\Cct}{\Zz[t,t^{-1}]}  
\newcommand{\LHZ}{\LH^{\Zz}}      


We consider a quotient algebra of the group algebra (over a suitable
commutative ring) of the group $\LBxn$.
The quotient algebra  
is 
\beq
\LHZ_\nn \; := \;
\faktor{\Cct\LBxn}{\R_n}
 \;\;  =\;\;   \langle \Xi_n \cup \Xi_n^{-} \; | \; \Q_n, \I_n, \R_n \rangle_{\Cct}
\eq
where $\R_n$ is the set of (algebra) relations:  
\begin{align}  \label{eq:R1} 
  \sigg\ii2 &= (1-t)\sig\ii + t
    \hspace{1in}\mbox{(i.e. $\; (\sig\ii-1)(\sig\ii+t)=0$)}
  \\
  \rr\ii\sig\ii &= -t\rr\ii + \sig\ii +t
    \hspace{1in}\mbox{(i.e. $\; (\rr\ii-1)(\sig\ii+t)=0$)}
\label{eq:R1i}
    \\ 
\sig\ii\rr\ii &= -\sig\ii + \rr\ii+1.
    \hspace{1in}\mbox{(i.e. $\; (\sig\ii-1)(\rr\ii+1)=0$)}
\label{eq:R1ii}
\end{align}
(NB we already have $(\rr\ii-1)(\rr\ii+1)=0$).

\ignore{{
\todopm{-Implicit in this seems to be that $t \in \Cc$. Is that what
  we have in mind? It does not exactly match what happened above...
  (this can all be sorted out, but it would be safer to settle it). }
\ppm{OK I'm going to try to do this now...}
}}

Observe that \eqref{eq:R1} yields an inverse for $\sig\ii$
(the inverse to $t$ is specifically needed),
so we have
the following:
\beq
\LHZ_n =   \;\; \langle \Xi_n  \; | \; \Q_n,  \R_n \rangle_{\Cct}
\eq
Observe then that the relations as such do not require an inverse to
$t$, so we could consider the variant algebra over $\Zz[t]$. 

For any field $K$ that is a $\Cct$ algebra
we then define the base change $\LH^K_n = K \otimes_{\Cct} \LHZ_n$
and,
for given $t_c \in \Cc$,
\[
\LH_n(t_c) = \LH_n = \LH^{\Cc}_n
\]
where $\Cc$ is a $\Zz[t]$-algebra by evaluating $t$ at
$t_c$  
(the choice of which we notationally suppress).
Note that there is no reason to suppose that this gives a flat
deformation (i.e. the same dimension) in all cases. 
(It will turn out that it does,
at least in low rank,
if we can localise at $t^2 -1$.
In particular, perhaps surprisingly, in the variant $t=0$ is
isomorphic to the generic case.)

\medskip\bigskip

\noindent
{\em Remark.} The relations \eqref{eq:R1} et seq
are suggested by  
\eqref{eq:loc} and 
the following calculations (on $\sig1$ and $\rho_1$ in $\LB3$,
noting that blocks work the same way for all generators):

\ignore{{
\begin{equation*}
\rrb{GB}(\sigg12)=
    \begin{pmatrix} 
     t^2-t+1 & t-t^2 & 0\\
     1-t & t & 0\\
     0 & 0 & 1\\
    \end{pmatrix} 
 = 
(1-t)\begin{pmatrix} 
     1-t & t & 0\\
     1 & 0 & 0\\
     0 & 0 & 1\\
     \end{pmatrix} 
+ 
t I_3,
\end{equation*}
}}

\begin{equation*}
\rrb{GB}(\sig1\rr1)=
    \begin{pmatrix} 
    t & 1-t & 0\\
    0 & 1 & 0 \\
    0 & 0 & 1 \\
    \end{pmatrix} 
=
-   \begin{pmatrix} 
     1-t & t & 0\\
     1 & 0 & 0\\
     0 & 0 & 1\\
    \end{pmatrix} 
+
    \begin{pmatrix} 
     0& 1 & 0\\
     1 & 0 & 0\\
     0 & 0 & 1\\
    \end{pmatrix}   
+ I_3
\end{equation*}


\begin{equation*}
\rrb{GB}(\rr1\sig1)=
    \begin{pmatrix} 
    1 & 0 & 0\\
    1-t & t & 0 \\
    0 & 0 & 1 
    \end{pmatrix}   
=
- t \begin{pmatrix} 
     0& 1 & 0\\
     1 & 0 & 0\\
     0 & 0 & 1
    \end{pmatrix}   
+
    \begin{pmatrix} 
     1-t & t & 0\\
     1 & 0 & 0\\
     0 & 0 & 1
    \end{pmatrix}   
+t I_3.
\end{equation*}


\subsection{Notable direct consequences of the relations: Finiteness}

\newcommand{\steq}[1]{\stackrel{#1}{=}}

Given a word in the generators, of form $\sig3\sig4\rr2$ say, by a
{\em translate} of it we mean the word obtained by shifting the indices
thus:
$\sig{3+i}\sig{4+i}\rr{2+i}$. 

With the $\Q$ and $\R$ relations we can derive the following ones,
together with the natural translates thereof
(here $\steq{*}$ uses \eqref{eq:pp1}; $\steq{\rr{}\rr{}\sig{}}$ uses
\eqref{eq:Qm0}; $\steq{\sig{}\rr{}}$ uses \eqref{eq:R1ii}, and so on): 
\ignore{{
\begin{itemize}
\item[(M1)]
  ${\rr\bb\sig\cc\rr\bb}
  =
  \rr2\sig1\rr2\rr1\rr1=\rr2\rr2\rr1\sig2\rr1=
  {\rr1\sig2\rr1}$,
\item[(M2)]
  $
{\sig2\rr1\rr2}
  =\rr1\rr1\sig2\rr1\rr2 =
  \rr1\rr2\sig1\rr2\rr2
  = 
  {\rr1\rr2\sig1}$
  \ppm{This is already in $\Q$! Don't need the middle bits.},
\end{itemize}

Using also the $\R$-relations we have the translates of:
}}
\begin{enumerate}
\item[(M1)]
  $\displaystyle\begin{aligned}[t]\mathbf{\sig2\rr1\sig2}
  &\steq{*} \sig2\rr2\rr2\rr1\sig2
  \steq{\rr{}\rr{}\sig{}} \sig2\rr2\sig1\rr2\rr1
  \steq{\sig{}\rr{}} -\sig2\sig1\rr2\rr1+\rr2\sig1\rr2\rr1+\sig1\rr2\rr1 \\
  &\steq{\rr{}\sig{}\sig{},\rr{}\rr{}\sig{}} -\rr1\sig2\sig1\rr1+\rr2\rr2\rr1\sig2+\sig1\rr2\rr1
  \steq{\sig{}\rr{}} \mathbf{\sig1\rr2\rr1+\rr1\sig2\sig1-\rr1\sig2\rr1},\end{aligned}$
\item[(M2)]
  \label{Rel:forbiddenn}
  $\displaystyle\begin{aligned}[t]\mathbf{\rr2\sig1\sig2}
  &\steq{*}\rr2\sig1\rr2\rr2\sig2
  \steq{\rr{}\sig{}} -t\rr2\sig1\rr2\rr2+\rr2\sig1\rr2\sig2+t\rr2\sig1\rr2
  \steq{*,\rr{}\sig{}\rr{}} -t\rr2\sig1+\rr1\sig2\rr1\sig2+t\rr1\sig2\rr1 \\
  &\steq{M1} -t\rr2\sig1+\rr1(\rr1\sig2\sig1-\rr1\sig2\rr1+\sig1\rr2\rr1)
  +t\rr1\sig2\rr1\\
  &\steq{*,\rr{}\sig{}} -t\rr2\sig1+\sig2\sig1-\sig2\rr1+(-t\rr1+\sig1+t)\rr2\rr1+t\rr1\sig2\rr1\\
  &= 
  {\sig1\rr2\rr1+t\rr1\sig2\rr1-t\rr1\rr2\rr1+\sig2\sig1
    -\sig2\rr1-t\rr2\sig1+t\rr2\rr1}.
  \end{aligned}$
\end{enumerate} 

\newcommand{\LHs}{\LH^{\langle}}

\begin{definition}
  For given $n$ and $m \leq n$ let $\LHs_m$ denote the subalgebra of $\LH_{n+1}$
generated by $\Xi_m$
(it is a quotient of $\LH_m$, as per the $\Psi$ map formalism in \S\ref{ss:arith1}).
\end{definition}

\begin{lemma} \label{lem:LXL}
For any $n$ 
  let $X_i$ be the vector subspace of $LH_{n}$ spanned by
  $\{ 1 ,\sig{i},\rr{i} \}$.
  Then $LH_{n+1} = \LHs_n X_n \LHs_n$. 
\end{lemma} 
\proof
It is enough to show that $ X_n \LHs_n X_n  $
lies in $\LHs_n X_n \LHs_n$.
We work by induction on $n$. The case $n=1$ is clear, since $LH_1 = \Cc$.
Assume true in case $n-1$ and consider case $n$. We have
$X_n \LHs_n X_n = X_n \LHs_{n-1} X_{n-1} \LHs_{n-1} X_n$ by assumption.
But $\LHs_{n-1}$ and $X_n$ commute so we have
$\LHs_{n-1} X_n X_{n-1} X_n \LHs_{n-1}$.
The inductive step follows from the relations $\Q$ and $\R$ and the relations
(M1,2) above.
\endproof

\begin{corollary}
$\LH_n$ is finite dimensional.  \qed
\end{corollary}

\begin{remark}\label{remark:other quotients}
  We may also 
  treat certain other quotients of $\C LB_n$.  For example,
eliminating either relations (\ref{eq:R1i}) or (\ref{eq:R1ii}) we still
obtain finite dimensional quotients.
In particular, if we only include
(\ref{eq:R1ii}) and not (\ref{eq:R1i}) then the analogous proof with $X_n$
replaced by
$\{1,\rr{n},\sig{n},\rr{n}\sig{n} \}$ proves finite dimensionality.
\end{remark}



\newcommand{\p}{\rho}
\newcommand{\s}{\sigma}
\newcommand{\tw}{t}  

\subsection{Refining the spanning set}
Can we express elements of $\LH_3$ as sums of length-2 words
(and hence eventually solve word problem)?
We have, for example, 
\beq \label{eq:ppp}
\p_1 \p_2 \p_1 =
-1 +\p_2
+\frac{(-\tw-1)}{(\tw-1)} (-\p_1 + \p_2 \p_1 -\p_1 \p_2)
+\frac{2}{(\tw-1)} (-\s_1 + \s_2 \p_1 -\p_1 \s_2)
\eq
But in general this is not easy. And another problem is that we do not
have immediately manifest relationships between different ranks
(such as inclusion) that would be useful. 
With this (and several related points) in mind it would be useful to have a tensor space
representation.
In what follows we address the construction of such a representation.



\section{Basic arithmetic with $\LH_n$} \label{ss:basic}
Here we briefly report some  
basic arithmetic in $\LH_n$ that gives
the clues we need for our local representation constructions below. 




\subsection{Fundamental tools, locality}

In what follows, $\Bcat$ denotes the {\em braid category}:
a strict monoidal category
with object monoid $(\Nn_0,+)$ generated by 1,
and $\Bcat(n,n) = B_n$,
$\Bcat(n,m)=\emptyset$ otherwise, and monoidal composition is
via side-by-side concatenation of suitable braid representatives
(see e.g. \cite[XI.4]{MacLane}).
Similarly $\Scat$ is the permutation category (of symmetric groups).
Let $\Hcat$ denote the ordinary Hecke category
--- again monoidal, but less obviously so \cite{Humphreys}. 
(We have not yet shown that $\LHcat$, the loop-Hecke category, is monoidal.)

Let $\LBcat$ denote the loop-braid
category
--- this is the strict monoidal category analogous to the braid
category
where the object monoid is $(\Nn,+)$, $\LBcat(n,n) = LB_n$,
$\LBcat(n,m)=\emptyset$ otherwise, and monoidal composition $\otimes$ is
side-by-side concatenation of loop-braids.

\bigskip

Suppose $\Ccat$ is a strict monoidal category with object monoid
$(\Nn_0 ,+)$ generated by 1 (for example, $\LBcat$).
Write $1_1$ for the unique element of $\Ccat(1,1)$
and
for $x \in \Ccat(n,n)$ define the  {\em translate}
\beq \label{eq:trans}
x^{(t)} = 1_1^{\otimes t} \otimes x \in \Ccat(n+t,n+t)
\eq
For $k$ a commutative ring, 
define translates of elements of $kLB_n$
(i.e. $k \LBcat(n,n)$), and $k S_n$ and so on,
by linear extension.

\bigskip

{\em Caveat:}
Note that it is a property of the geometric topological construction
of loop braids that the  
composition $\otimes$ in $\LBcat$  makes manifest sense. It requires
that side-by-side concatenation of rank $n$ with rank $m$ passes to
$n+m$. This is clear by construction.
\ignore{{
\footnote{Maybe say a bit more now?
  Observe that every class=loop-braid has a `narrow' representative,
  and we can use these.
}
}}
But in groups/algebras defined
by generators and relations it would not be intrinsically clear.
For example, how do we know that the subalgebra of $LH_n$ generated
{\em in} $LH_n$ by the elements
$p_i, s_i$, $i=1,2,...,n-2$ is isomorphic to $LH_{n-1}$?
(Some of our notation  
requires care
at this point since it may lead us to
take isomorphism for granted!)




\newcommand{\kalg}[1]{\langle #1 \rangle_k}
\newcommand{\XS}{X}  

\subsection{The $\Psi$ maps}
\label{ss:arith1}

Let
$A=\kalg{ \XS | R}$  
be an algebra
presented with generators $\XS$ and
relations $R$.
Then there is a homomorphism from the free algebra generated by any
subset $\XS_1$ of $\XS$ to $A$, taking $s \in \XS_1$ to its image in $A$.
This factors through the quotient by any relations, $R_1$ say, expressed only in
$\XS_1$.
We may consider it as a homomorphism from this quotient. But of course
the kernel may be bigger --- relations induced indirectly by the
relations in $R$.
A $\Psi$ map
is 
such a homomorphism:
\[
\kalg{ \XS_1 | R_1 }
\stackrel{\Psi}{\twoheadrightarrow}  \kalg{ \XS_1 | R }
\hookrightarrow  
\kalg{ \XS | R }
\]

Note that arithmetic properties such as idempotency, orthogonality and vanishing
are preserved under $\Psi$ maps.
Thus for example a decomposition of 1 into orthogonal idempotents in
$kS_n$ passes to such a decomposition in $LH_n$
(see \eqref{eq:PsiSLH}).
However conditions such as primitivity,
inequality and even non-zero-ness are not preserved in general.


Note that there is a natural (not generally isomorphic) image of
\beq 
kS_n \; \cong \; 
k\langle \rrp{1},...,\rrp\ii,...,\rrp{n-1} \;|\;
\A_n(\rrp{})    
, {\rrp\ii}{\rrp\ii}=1   
\rangle
\eq
in $LH_n$ obtained by the 
map of generators $\rrp{i} \mapsto \rr{i}$.
Let us call it $LH_n^{\rr{}}$. Thus
\beq \label{eq:PsiSLH}
k S_n \stackrel{\Psi}{\twoheadrightarrow} LH_n^{\rr{}} \hookrightarrow LH_n
\eq
Similarly
$H_n = \langle \sigT{1},...,\sigT\ii,...,\sigT{n-1} | \;\A_n(\sigT{}),... \rangle_k$
has image $LH_n^{\s}$
under $\sigT{i} \mapsto \sig{i}$:
\beq
H_n \twoheadrightarrow LH_n^\s \hookrightarrow LH_n
\eq

\bigskip

\ignore{{
By these inclusions we may build elements of $LH_n$ that will have
nice properties with respect to the actions of one subset of
generators or the other, derived from constructions for the
algebras on the left.

\bigskip
}}



Let us consider the image of a primitive idempotent  decomposition in $kS_n$:
\[
1 = \sum_{\lambda\in \Lambda_n} \sum_{i=1}^{d_\lambda} \e_\lambda^i
\]
under $\Psi : kS_n \rightarrow LH_n$.
Here $\Lambda_n$ denotes the set of integer partitions of $n$,
and $d_\lambda$ is the dimension of the $S_n$ irrep.
See \S\ref{ss:SnIds} for explicit constructions.
We will also write $(\Lambda, \subseteq)$ for the poset of all integer
partitions ordered by the usual inclusion as a Young diagram.


\begin{proposition} \label{pr:key22}
  Let $k$ be the field of fractions of $\Zz[t,t^{-1}]$.
  \\
(I) The image $\Psi(\e_\lambda^i)$ in $\LH^k_n$ of every idempotent with
  $(2,2) \subseteq \lambda \in \Lambda_n$
  is zero. 
\\
  (II) On the other hand
  all other $\lambda\in \Lambda_n$, i.e. all hook shapes, give
  non-zero image.
\end{proposition}
\proof
(I) 
Note that $\e_\mu^1$ with $\mu\in\Lambda_m$ is defined in $k S_n$ for
$n \geq m$ by $S_m \hookrightarrow S_n$.
It is shown for example in \cite{Rittenberg92} that if the relation
$\e_\mu^1 =0$ is imposed in a quotient of $k S_n$ then
$\e_\nu^i =0$ holds for $\mu\subseteq\nu\in \Lambda_n$
(a proof uses $S_{n-1}\hookrightarrow S_n$ restriction rules,
{from which we see that $\e_{\mu}^1$ is expressible as
  a sum of orthogonal such idempotents}).
%
Consider $\e^1_{(2,2)}$ (i.e. with $(2,2) \in \Lambda_4$) which may be expressed as
\[
\e^1_{(2,2)} \;= \;\;\; \raisebox{-.21in}{\includegraphics[width=1.52cm]{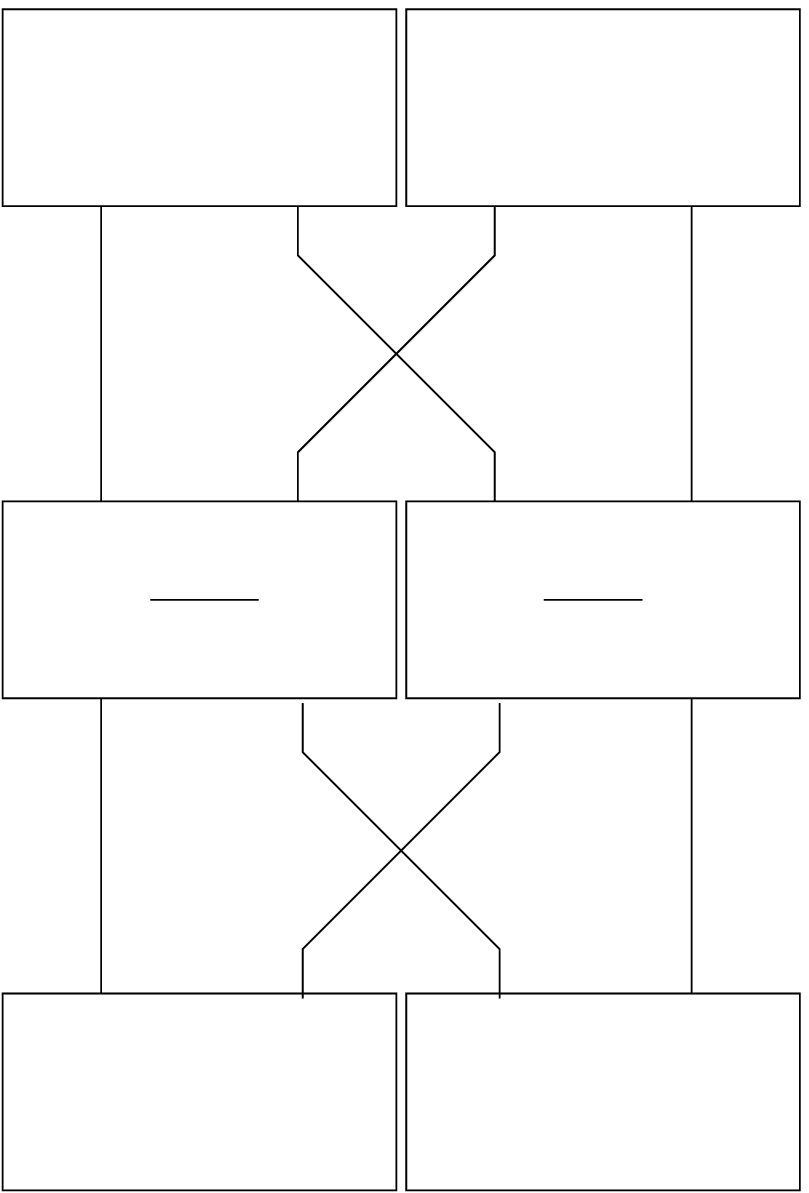} }
\; \propto \; 
( p_1 +1)(p_3+1) p_2 (p_1 -1)(p_3 -1) p_2 ( p_1 +1)(p_3+1)
\]
(using notation and a choice from (\ref{eq:choicee})).
By a direct calculation
in $LH_4$
\beq \label{eq:e220}
\Psi(( p_1 +1)(p_3+1) p_2 (p_1 -1)(p_3 -1)) = 0
\eq
(N.B. we know no elegant way to do this calculation;
the result holds also for generic $t$, but {\em not} for $t=1$).
\ignore{{
\footnote{{
Let us drop the $\Psi$ where no ambiguity arises. For example in $\LH^k_n$ from Magma:
\[
( p_1 +1)(p_3+1) x (p_1 -1)(p_3 -1) = 0
\]
for $x=1, \; p_2, \; s_2, \; p_2 s_2 , \; p_2 s_1 p_2, \; s_2 p_1 s_2$,
but {\em not} for $x=s_3 s_2, \; s_3 s_1, \; s_2 p_1 p_2$
(although all these might be linearly dependent).
}}
}}

\noindent
(II)
This can be verified by evaluation as non-zero in a suitable representation.
(For simplicity it is sufficient  to work in  the `SP quotient' that
we give in  Prop.\ref{th:ESSP} below, working with Kronecker products.
We will omit the explicit calculation.)
\qed


\medskip

With identity \eqref{eq:e220} in mind, recall  that in \cite{Rittenberg92}
local representations of ordinary Hecke
(and hence $S_n$) 
with this property were
constructed from spin chains.
In \S\ref{ss:locreps} we will combine this with Burau
and thus find 
the representations of loop-Hecke that we need here.

\medskip

By Prop.\ref{pr:key22} we have a decomposition of 1 in $\LH_n$ according to hook
partitions:
\beq \label{eq:d1} 
1 \;=\; 
\sum_{i=0}^{n-1} \; \sum_{j=1}^{d_{(n-i,1^i)}} \Psi( \e^j_{(n-i,1^i)} )
\eq
(NB $j$ varies over idempotents that are equivalent in the sense that
they induce isomorphic modules --- it will be sufficient to focus on $j=1$).
\ignore{{
Given the idempotent inclusion property \ref{prop:ii} in the SP section
below, the number of classes of these agrees with the number of simple
modules of $\LH_n$.
We can use these idempotents to finish the analysis.
We need to consider the action of {\em other} generators from $\LH_n$,
as follows.
}}

(Left) multiplying by $A=\LH_n$ we thus have a decomposition of the algebra
\[
A \;\; \cong \;\;  
\oplus_{i=0}^{n-1} \oplus_j A \Psi( \e^j_{(n-i,1^i)} )
\]
as a left-module for itself, into projective summands.

We have not yet shown that these summands are indecomposable.
\ignore{{
So for now we have
\[
A \Psi( \e^j_{(n-i+1,i-1)} ) = P \oplus Q \oplus ...
\]
where $P,Q,...$ are indecomposable projective modules.
(NB We claim that this module is in fact indecomposable, so there is
only one summand.
This is the same as saying that $\e^j_{(n-i+1,i-1)} $ remains
primitive in $\LH_n$ as it is in $\Cc S_n$.)
}}
But consider for a moment the action of $\LH_n$ on the image under
$\Psi$ of
$$
Y^n_{\pm} =  \sum_{g \in S_n} (\pm 1)^{len(g)} g
$$
in
$\LH_n$
(we write $Y^n_+$ for unnormalised $\e^1_{(n)}$ and $Y^n_-$ for $\e^1_{(1^n)}$
--- again see \S\ref{ss:SnIds} for a review).
By abuse of notation we will write $Y^n_{\pm}$ also for the image. 
By (\ref{eq:R1ii}) and
the classical identities $ Y^{a(1)}_{\pm} Y^n_{\pm} = a! Y^n_{\pm}$ 
(recall $ Y^{a(1)}_{\pm} $ means $ Y^{a}_{\pm}$ with indices shifted
by +1, 
see (\ref{eq:Yid1}) {\em et seq})
we have
\beq \label{eq:sigY}
\sig\ii Y^n_{+} = Y^n_+ ,  \hspace{1.9in}
Y^n_{-} \; \sig\ii  \;\; =\;  -t Y^n_{-}
\eq
It follows that $Y^n_+$ spans a 1-d left ideal in $LH_n$.
If we work over a field containing the rationals then it is
normalisable as an idempotent, and so we have an indecomposable
projective left module
\[
P_{(n)} = \; \LH_n Y^n_+ \; = \; \LH_n \e^1_{(n)}  \;
       = \; k \e^1_{(n)} .
\]


\section{On local representations} \label{ss:locreps}





Here $\Mat$ is the monoidal category of matrices over a given
commutative ring (and $\Mat_k$ the case over commutative ring $k$),
with
object monoid $(\Nn, \times)$ and
tensor product on morphisms given by  
a Kronecker product
(NB there is a convention choice in defining the Kronecker product).
We often focus on the monoidal subcategory $\Mat^m$ generated  by a
single object $m \in \Nn$ --- usually $m=2$.
Then the object monoid  $(2^{\Nn}, \times)$ 
becomes $(\Nn,+)$ in the natural way.



\newcommand{\rb}[2]{\raisebox{#1pt}{#2}}
\newcommand{\BOX}{\rb{-4}{\includegraphics[width=1cm]{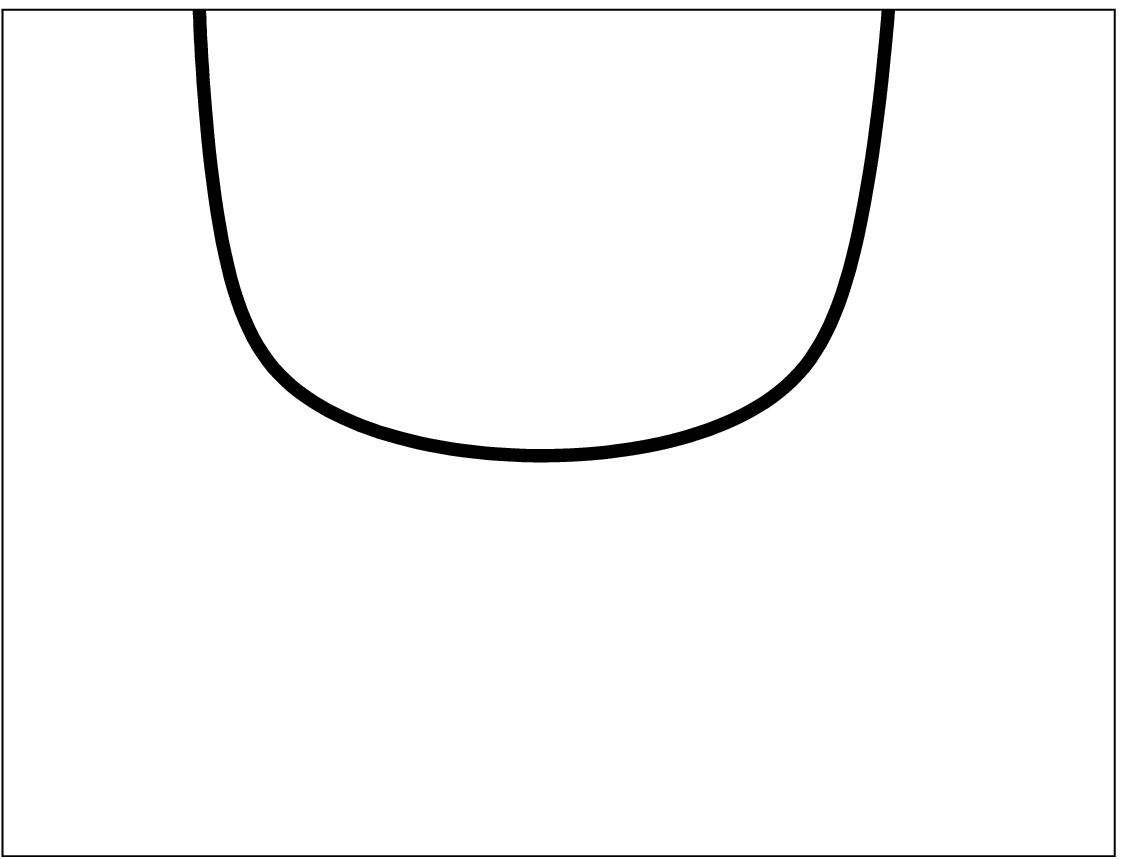}}}
\newcommand{\DUAL}{\rb{-4}{\includegraphics[width=1cm]{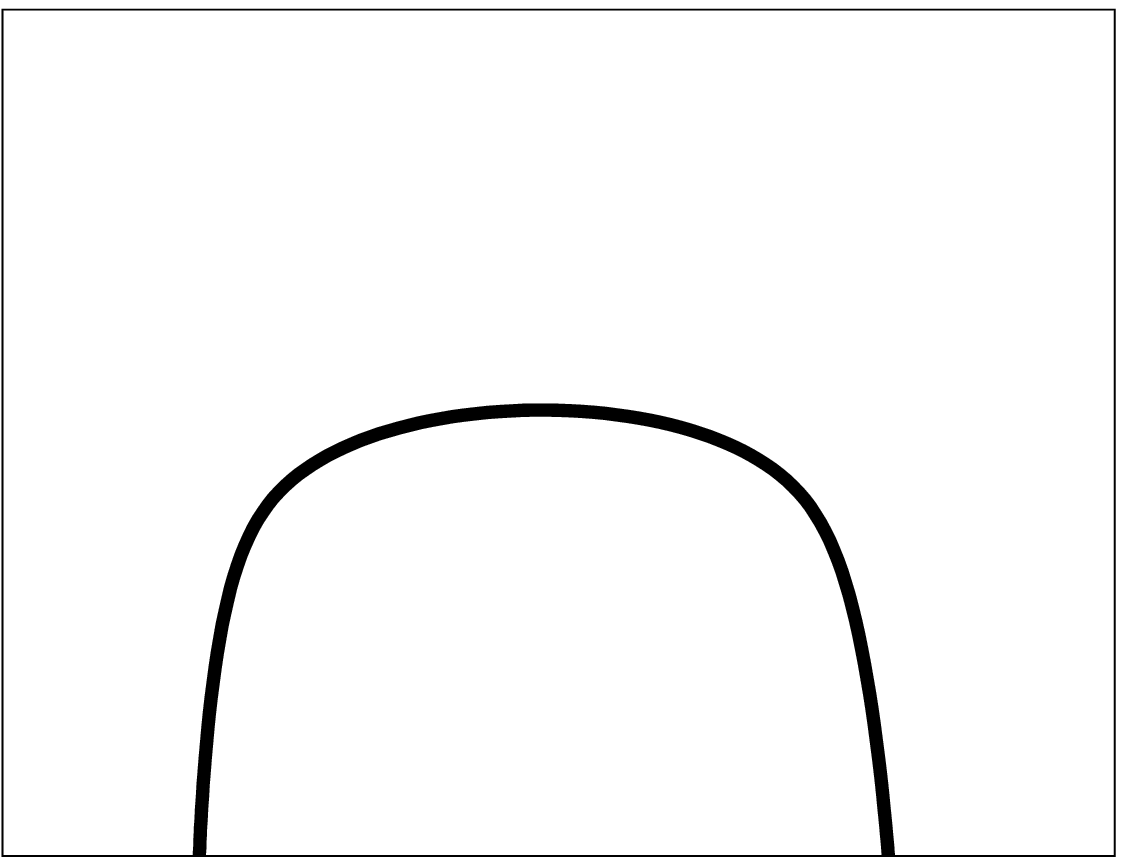}}}
\newcommand{\ud}{{\mathsf u}}
\newcommand{\Tcat}{{\mathsf T}}
\newcommand{\ttt}{\tau}
\newcommand{\U}{{\mathsf U}}
\newcommand{\G}{{\mathsf G}}
\newcommand{\GU}{\G_{\U}}
\newcommand{\ups}{\Upsilon}

In the study of ordinary Hecke algebras (and particularly
quantum-group-controlled quotients like Temperley--Lieb) a very useful
tool is the beautiful set of local tensor space representations
generalising those arising from XXZ spin chains and Schur--Weyl
duality.
For example we have the following.

Consider the TL diagram category $\Tcat$
with object monoid $(\Nn,+)$
$k$-linear-monoidally generated by the
morphisms represented by diagrams
$\ud = \BOX \in \Tcat(2,0)$ and $\ud^* = \DUAL$.
This has a TQFT
$\FF_2$
given by
$\BOX \mapsto (0, \ttt , \ttt^{-1} , 0)$
(the target category is $\Mat$)
and  taking $*$ to transpose.
Of course for $1_1 \in \Tcat(1,1)$ we have $\FF_2( 1_1) = I_2$.

To pass to our present topic we note
that $1_1 \otimes 1_1 = 1_2$ and that 
the Yang--Baxter construction
$\sigma_1 \mapsto 1_2-\ttt^2 \ud^* \ud$
gives
\beq \label{eq:YB}
\sigma_1 \; \mapsto \;\; 
\FF_2( 1_2 )
-\ttt^2 \mat{ccccc} 0 \\ &\ttt^2 &1 \\ &1 &\ttt^{-2} \\ &&& 0 \tam
\;\; = \mat{ccccc} 1 \\ &1-\ttt^4 &-\ttt^2 \\ &-\ttt^2 &0 \\ &&& 1 \tam
\eq
thus a representation of
the braid category $\Bcat$
(note that eigenvalues are 1 and $-\ttt^4$ so 
$\ttt^4$ here passes to $t$ in our parameterisation for loop-Hecke).
But also note that $\ud ,\ud^* $ can be used for a Markov trace.
And also for idempotent localisation functors:
let $\U = \ud^* \ud$, $\U_1 = \U \otimes 1_{n-2}$,
and $T_n = \Tcat(n,n)$ regarded as a $k$-algebra;
then we have the algebra isomorphism
$\U_1 T_n \U_1 \cong T_{n-2} $.
This naturally gives a  category embedding $\GU$ of $T_{n-2}$-mod in $T_n$-mod.
Recall that irreps are naturally indexed by partitions of $n$ into at
most two parts: $\lambda = (n-m,m)$, or equivalently 
(for given $n$)
by `charge' $\lambda_1 - \lambda_2 =n-2m$,
thus by $\ups_n = \{ n,n-2, n-4,...,0/1 \}$ (depending on $n$ is odd or even).
This latter labelling scheme is stable under the embedding.
That is, indecomposable projective modules are mapped by $\GU$
according to $\ups_{n-2} \hookrightarrow \ups_n$. 
%
\ignore{{
\footnote{
Note for later that in case $\ttt^4 = -1$ (FIX ME!)
$$
\sigma_1 \mapsto 
\mat{ccccc} 1 &&& 1 \\ &1-\ttt^4 &-\ttt^2 \\ &-\ttt^2 &0 \\ &&& 1 \tam
$$
also
almost
gives a solution to the YB relations. That is ... 
}
}}

\subsection{Charge conservation}

\newcommand{\bbb}[1]{\varepsilon_{#1}}  

Another useful property of $\FF_2$ is `charge conservation'.
We may label the row/colum index for object 2 in $\Mat$ by
$\{\bbb{1},\bbb{2} \}$ or $\{ +,- \}$. Then $2\otimes 2$ has index set
$\{ \bbb{1}\otimes \bbb{1},
\bbb{2}\otimes \bbb{1}, \bbb{1}\otimes \bbb{2}, \bbb{2}\otimes \bbb{2} \}$
(which we may abbreviate to $\{ 11,21,12,22 \}$)
and so on.
The `charge' $ch$ of an index is $ch = \#1-\#2$.
Note from \eqref{eq:YB} that $\FF_2$ does not mix between different
charges (hence charge conservation).

\newcommand{\UU}{{\mathcal U}}  %
\newcommand{\UUT}{{\mathcal T}}  %
\newcommand{\eb}{\beta}  

For a functor with the charge conservation property the representation
of $B_n$ (say) obtained has a direct sum decomposition according to
charge, with
`Young blocks'
$\eb_i$ of charge $i=n,n-2,...,-n$.
The dimensions of the blocks are given by Pascal's triangle.
It will be convenient to express this with the semiinfinite
Toeplitz matrices
$\UU$ and $\UUT$:
\[
\UU = \mat{cccccccccc} 1&1 \\ &1&1 \\ &&1&1 \\ &&&1&1 \\ &&&& \ddots \tam
, \hspace{.4602cm}
\UU^2 = \mat{cccccccccc} 1&2&1 \\ &1&2&1 \\ &&1&2&1 \\ &&&1&2&1 \\ &&&& \ddots \tam
, \hspace{.471cm}
\UUT = \mat{cccccccccccc} 0 & 1 & \\ 1&0&1 \\ &1&0&1 \\ &&1&0&1 \\ &&&& \ddots \tam
\]
and semiinfinite vectors $v_1 = (1,0,0,0,...)$,
$v_2 = (0,1,0,0,...)$, ....
Thus $v_1 \UU^n$
(respectively $v_{n+1} \UUT^{n}$)
gives
the numbers in the $n+1$-th row of Pascal (followed by a tail of
zeros).
(The two different formulations correspond to two different
thermodynamic limits
---  $\UUT$ corresponds to the $\ups_{n-2} \hookrightarrow \ups_n$ limit
--- see later.)
Then
\beq \label{eq:Pascn}
dim(\eb_i)  \; = \;  \left( v_1 \UU^n \right)_{(n-i+2)/2}
       \; = \; \left( v_{n+1} \UUT^{n} \right)_{n-i+1}
\eq
In the case of $\FF_2$ these blocks are not linearly irreducible in
general (the generic irreducible dimensions are given by $v_1 \UUT^n$).
But they still provide a useful framework. We return to this
later.

With this construction and Prop.\ref{pr:GBx} in mind, it is natural to ask if we
can make a local version of generalised Burau.
(Folklore is that this cannot work, and directly speaking it does not.
But we now have some more clues at our disposal.)



\subsection{Representations of $\Bcat$}  

Now we have in mind Prop.\ref{pr:key22}; and 
brute force calculations in low rank showing (see \S\ref{ss:lowrank})
that
$\LH_n$ is
non-semisimple but has irreducible representations with dimensions given by Pascal's
triangle.
This is  
reminiscent of 
Rittenberg's analysis of
Deguchi {\em et al}'s
quantum spin chains over Lie superalgebras
\cite{deguchi_89,deguchi1989,Deguchi_1990,kauffman1991,Rittenberg92,Deguchi92}
(it is also reminiscent of work of Saleur on `type-B' braids
\cite{Martin-Saleur:M3a}
--- but for this
cf. e.g. \cite{Bullivant-et-al:2020}).  
Inspired by    this     
and the Burau representation
(and cf. \cite{Damiani2018} and references therein) we  
proceed as follows.
%
%
%
Define
\beq
M_t(\sigma)  = \mat{cccc} 1 \\ & 1-t & t  \\ & 1 & 0  \\ &&& 1 \tam
,\hspace{.61in} 
M'_t(\sigma)  = \mat{cccc} 1 \\ & 1-t & t  \\ & 1 & 0  \\ &&& -t \tam
\eq
as in \cite{deguchi_89,kauffman1991}.
\ignore{{
\begin{proposition}
  \ppm{{\rm (\cite{deguchi_89,kauffman1991})}
  }
  }}
Fix a commutative ring $k$,  $\tau\in k^\times$,  and $t=\tau^4$.
Observe that there 
is a monoidal functor $\FF_M$ from the
Braid category $\Bcat$ 
to $\Vect$ (or at least $\Mat$) given by
object 1 mapping to $V= \Cc\{ e_1,e_2 \}$
(i.e. to 2 in $\Mat_{\Zz[t]}$)
and
the positive braid $\sigma$  in $\Bcat(2,2)$ mapping to
$ M_t(\sigma) $.
The conjugation of this matrix to $\FF_2 (\sigma)$ lifts to a natural
isomorphism of functors.
Another natural isomorphism class of charge conserving functors
has representative 
functor  $\FF_{M'}$ given by $ M'_t(\sigma) $.
(According to  Deguchi {\em et al}'s  scheme this  
is the (1,1)-super class,
cf. for example \cite{deguchi_89,kauffman1991}. 
But note that in extending to $\LBcat$ below, isomorphism will not be preserved,
so we are focussing on the specific representative.) 
In fact
some elementary analysis shows that these two classes are all
of this form 
that factor through Hecke (apart
from the trivial class). 


Let us formulate this in language that will be useful later. First
note that (like any invertible matrix) $M_t(\sigma)$ and
$M'_t(\sigma)$
extend to monoidal functors from the free monoidal category generated
by $\sigma$ to $\Mat$.
Thus in particular
$M'_t(\sigma\otimes 1_1 ) =  M'_t(\sigma) \otimes Id_2  \in Mat(2^3,2^3)$.
Given the form of the construction, proof of the above
factorings through $\Bcat$ 
follows from a
direct verification of the braid relation in each case.
\ignore{{
\ppm{IM about TO DELETE THE NEXT BIT... RESUME AT `Indeed...' below.}
\proof
\ppm{
For well-definedness see e.g. Deguchi et al \cite{deguchi_89,deguchi1989,akutsu87}
or Kauffman--Saleur \cite{kauffman1991}.
However, since we will need to generalise later, we 
recall the following points.}
There is a proof
in each case  by  direct calculations
and properties of the construction.
Let $1_n$ denote the identity morphism in $\Bcat(n,n)$. Then 
\[
  M_t( (\sigma\otimes 1_1 )(1_1\otimes\sigma)(\sigma\otimes 1_1 ) )
  = M_t( (1_1\otimes\sigma)  (\sigma\otimes 1_1 )(1_1\otimes\sigma)  )
\]
with other relations then 
holding by construction.
Invertibility of the image of $\sigma$ is clear.
$\;$ 
The $M'$ case is similar.
The proof of completeness is elementary analysis.
\qed

\medskip


Indeed
}}
    More interestingly we have,
    again by direct calculation,
    the  stronger result:
\beq \label{eq:tts}
  M'_t (\sigma\otimes 1_1 ) M'_t (1_1\otimes\sigma) M'_s (\sigma\otimes 1_1 ) 
  = M'_s (1_1\otimes\sigma) M'_t (\sigma\otimes 1_1 ) M'_t (1_1\otimes\sigma)  
\eq
while the $tss$ version of this  identity does {\em not} hold
(unless we force $s=1$, or $s=t$) 
(NB care must be taken with conventions here). 

\ignore{{
\medskip

Recall
that the $M_t$ case is a smooth (not necessarily flat)
    deformation of the $t=1$ case, which is the classical $N=2$ tensor
    space representation 
    of $\Scat$ --- a special case of our $\FF_2$ functor in \eqref{eq:YB} above. 
  For $n=3$ the $M_t$ case 
  has 2 copies of the trivial rep and 2 copies
  of the Burau rep (thus dim: $2^3 = 1+3+3+1$).
  Thus at $t=1$ it is a direct sum of
  four copies of the trivial and two copies of the 2d  
  $S_n$ irrep.
  For $n=4$ we have 16=1+4+6+4+1, so as well as two trivial
  and two Burau (4=3+1), now also a $6$ (with $6=2+3+1$ when $t=1$)
  --- see e.g. \cite[\S9.6]{Martin:M1}.
  The $M_t'$ case has the same block decomposition but different
  irreducible decomposition \cite{Rittenberg92,Deguchi92}.

}}






\subsection{Extending to $\LBcat$}

Recall we introduced the loop-braid category $\LBcat$.
We write $\sigma \in \LBcat(2,2)$ for the positive braid exchange and $\rho\in
\LBcat(2,2)$ for the symmetric exchange.

Formally 
extending with  elementary transpositions (cf. $\varrho_{GB}$), 
the $\FF_M$ construction  fails to satisfy the mixed braid relation
\eqref{eq:Qm}.
\ignore{{
The reason for the experiment is simply to bring in Pascal
combinatorics in a naive way --- since the $(n-m,m)$ `fixed charge'
blocks have it
(for example $M$ itself has block shape 1+2+1).
Of course the failure is, in this linear block sense,
not complete --- relaxing the monoidal structure,
the $(n-1,1)$ subspace is exactly generalised Burau.
Can we adjust this to fix the failure
(and the braid rank shift --- note that $M$ is `rank 2' not rank 3 cf. (\ref{eq:Pascal}))
somehow?! Yes:
}}
However the 
functor $\FF_{M'}$ fairs better:


\begin{theorem} \label{th:ESSP}
  (i) The $\sigma\mapsto M'_t(\sigma)$
  construction extended using the super transposition
$\rho \mapsto M'_1(\sigma)$
  gives a monoidal functor
$F^e_{M'}$  from the loop  
Braid category 
$\LBcat$ to $\Mat$.
\\
(ii) $F^e_{M'}$ factors through $\LHcat$. 
\end{theorem}
\proof
The proof is a linear algebra calculation similar to the $\Bcat$ cases above,
using Kronecker product identities;
but also using
the appropriate special case of \eqref{eq:tts}
for \eqref{eq:Qm}.
\qed

\newcommand{\SSP}{SP} 

\begin{definition}
Fix a field $k$ and $t \in k$.
Then the $k$-algebra 
$\SPe_n = k LB_n / \ann F^e_{M'}$. 
\end{definition}


We conjecture that the extended super  
representation,
which we call Burau--Rittenberg, or 
`\SSP' rep for short,
is faithful
on $\LH$ unless $t^2 =1 $ (see later).

\begin{proposition}
  Fix a field $k$ and $t\in k$, $t \neq 1$.
  Let $\chi_i = \frac{\sig{i}-\rr{i}}{1-t}$.
  Then \\ (a) $\chi_i$ and $\rr{i}$
  ($i=1,2,...,n-1$)
  are alternative generators of
  $\SPe_n$; and
\\  (b) The $k$-algebra isomorphism class of $\SPe_n$ is independent of $t$.
\end{proposition}

\proof
(a) Elementary. (b) The images of the alternative generators in the defining
representation are independent of $t$.
 \qed

\subsection{Towards linear structure of SP}

Let us work out the linear structure of SP. (I.e. its Artin--Wedderburn linear
representation theory over $\Cc$: simple modules, projective modules
and so on.
See \S\ref{ss:repthy} for a review.)


\begin{proposition} \label{pr:SP2}
Suppose $t \neq 1 \in k$. 
  Let $\chi  = \frac{ \sigma_{} - \rho_{} }{1-t}$ and
$\chi_1  = \frac{ \sigma_1 - \rho_1 }{1-t} \in \SPe_n$.
Then
\beq \label{eq:chiSP}
\chi_1 \SPe_n \chi_1 \cong \SPe_{n-1}
\eq
and 
\beq \label{eq:SPchi}
\SPe_n \;/\; \SPe_n \;\chi_1\; \SPe_n \;\; \cong \; k .
\eq
(II) 
In particular the map
$f_\chi : \SPe_{n-1} \rightarrow  \chi_1 \SPe_n \chi_1$
given by $w \mapsto \chi_1 w^{(1)} \chi_1$
(recall the translation notation from (\ref{eq:trans}))
is an algebra isomorphism.
\end{proposition}
\proof
Let us write simply $F=F_n$ for the defining representation $F^e_{M'}$ of
$\SPe_n$. We write $\{1,2\}^n$ for the basis
(i.e. we write simply symbols $1,2$ for $e_1, e_2$;
and the word $112$ for $e_1 \otimes e_1 \otimes e_2$ and so on). 
Our convention for ordering the basis is given by 11,21,12,22.
First observe that the image in $F$ is (here with $n=3$)
\beq
(\chi\otimes 1_2)  \mapsto
\ignore{{
\mat{ccccccccc}
0 \\
 &\fb{1-t}&t-1 \\
 & &0 \\
 & & & \fb{1-t}& \\
 & & &    &0 \\
 & & & &&\fb{1-t}  &{t-1} \\
 & & &    & & &0 \\
 & & &    & & & &\fb{1-t}
\tam
}}
\mat{ccccccccc}
0 \\
 &\fb{1}&-1 \\
 & &0 \\
 & & & \fb{1} \\
\tam
\otimes 1_2 =
\mat{ccccccccc}
0 \\
 &\fb{1}&-1 \\
 & &0 \\
 & & & \fb{1}& \\
 & & &    &0 \\
 & & & &&\fb{1}  &{-1} \\
 & & &    & & &0 \\
 & & &    & & & &\fb{1}
\tam
\eq
Note that the basis change conjugating by
\beq \label{eq:basch}
\mat{ccccc} 1 \\ & 1&1 \\ &&1 \\ &&&1 \tam \otimes 1_2
\eq
(again this is the example with $n=3$)
brings this into diagonal form, projecting onto the $2\{1,2\}^{n-1}$
subspace
(the subspace of $V^n = \Cc \{ 1,2 \}^n$ spanned by basis elements
of form $2w$ with $w \in \{ 1,2 \}^{n-1}$,
i.e. of form $e_2 \otimes ...$).
That is:
$\chi \mapsto \mat{cc} 0 \\ &1\tam \otimes 1_2$.
Furthermore, 
%
%
\beq
(\chi\otimes 1_2) (1_2 \otimes \sigma) (\chi\otimes 1_2 ) \mapsto
\mat{ccccccccc}
0 \\
 &\fb{1}&-1 \\
 & &0 \\
 & & & \fb{1-t}&0&\fb{t}&-t \\
 & & &    &0 \\
 & & & \fb{1}  &0&\fb{0} \\
 & & &    & & &0 \\
 & & &    & & & &\fb{-t}
\tam
\eq
Note that after the (\ref{eq:basch}) basis change this decomposes as a
sum of several copies of the 0 module together with the submodule $N$ with
basis $2\{1,2\}^{n-1}$. Then the map
from $N$ to $\{ 1,2 \}^{n-1}$ given by 
$2w \mapsto w$ gives
\[
\chi_1 \sigma_1 \chi_1 = -t.\chi_1 \mapsto F_{n-1}(-t.1)
\]
and $\chi_1 \sigma_2 \chi_1  \mapsto F_{n-1}(\sigma_{1})$
and $\chi_1 \rho_2 \chi_1  \mapsto F_{n-1}(\rho_{1})$.
Also note that $\chi_1$ commutes with $\sigma_i$ for $i>2$ so we have
\beq \label{eq:fubar}
\chi_1 \sigma_i \chi_1 = \chi_1 \sigma_i  \mapsto F_{n-1}(\sigma_{i-1}) ,  \hspace{1cm} 
\chi_1 \rho_i \chi_1  \mapsto F_{n-1}(\rho_{i-1})  \hspace{1cm} i>2 .
\eq
Thus the images of the generators under $w \mapsto \chi_1 w \chi_1$
are the generators of $\SPe_{n-1}$, establishing (\ref{eq:chiSP})
on generators.
To show that the images of the generators span we proceed as follows.
From Lemma~\ref{lem:LXL} and the
(sufficient)
symmetry of the relations under
$i \mapsto n-i$ on indices, writing $\LH_n = L_n$ for short, we have
\[
L_{n+1} = L_n X_n L_n  = L_n^{(1)} X_1 L_n^{(1)}
= L_{n-1}^{(2)} X_2 L_{n-1}^{(2)}  X_1 L_{n-1}^{(2)} X_2 L_{n-1}^{(2)}
\] \[
= L_{n-1}^{(2)} X_2 L_{n-2}^{(3)}  X_3 L_{n-2}^{(3)} X_1 X_2 L_{n-1}^{(2)} 
= L_{n-1}^{(2)} X_2   X_3  X_1 X_2 L_{n-1}^{(2)} 
\]
Thus
\[
\chi_1 L_{n+1} \chi_1 =
\chi_1 L_{n-1}^{(2)} X_2   X_3  X_1 X_2 L_{n-1}^{(2)} \chi_1
=
L_{n-1}^{(2)} \chi_1 X_2   X_3  X_1 X_2 \chi_1 L_{n-1}^{(2)} 
\]
We can show by direct calculations that
$\chi_1 X_2   X_3  X_1 X_2 \chi_1$ lies in the algebra generated by
the
images of the generators.
(We can do this even in $\LH_4$. The result then holds in $SP_4$ since it
is a quotient; and then in $SP_n$ by construction.
--- Note however that we have not shown that it holds in $\LH_n$.)
Also $L_{n-1}^{(2)} \chi_1  $ evidently lies in the algebra generated
by the images of the generators, by commutation,
so we are done. 

Finally (\ref{eq:SPchi}) follows on noting that
the quotient corresponds to imposing $\chi_1=0$, i.e. $\sigma_1 =
\rho_1$.
Noting that $t \neq 1$, this gives $\sigma_i =1$. 
\\
(II) Note that $f_\chi$ inverses the map from (\ref{eq:fubar}) above.
\qed

\subsection{Aside on linear/Artinian representation theory} \label{ss:repthy}

Since this paper bridges between topology and linear representation
theory it is perhaps appropriate to say a few words on the bridge.
While topology focuses on topological invariants, linear rep theory
is concerned with invariants such as the spectrum of linear
operators (and the generalised `spectrum' of algebras of linear operators).
The former is thus of interest for topological quantum
field theories, and the latter for usual quantum field theories
(where notions such as mass are defined).
In this section we recall a few key points of linear/Artinian rep
theory that are useful for us.
(So of course it can be skipped if you are not interested in this
aspect, or are already familiar.) 

Recall that every finite dimensional algebra over an algebraically
closed field is Morita equivalent to a basic algebra
(see e.g. \cite{NesbittScott,Jacobson,Benson95}).
This allows us to track separately the combinatorial and homological
data of an algebra.

Let $A$ be a finite dimensional algebra over an algebraically closed
field $k$ (cf. e.g. \cite{Benson95}).
Let $J(A)$ denote the radical. 
Let $L= \{ L_1, ... L_r \}$ be an ordered set of the isomorphism classes of simple $A$-modules,
with projective covers $P_i = Ae_i$
(i.e. the $e_i$s are a set of primitive idempotents).
Given an $A$-module $M$ let $Rad(M)$ denote the intersection of the maximal
proper submodules.
Now suppose $A$ is basic.
Recall that $Ext^1_A(L_i,L_j)$ codifies the non-split extensions
between these modules --- i.e. the `atomic' components of
non-semisimplicity. 
The corresponding `Ext-matrix' $E_L(A)$ is given by
\[
( E_L(A) )_{ij} = \; dim_k Ext^1_A(L_i, L_j)
\]
or equivalently
\[
dim_k Ext^1_A(L_i, L_j) = \;
dim_k ( Hom_A(P_j,Rad(P_i))/Hom_A(P_j,Rad^2(P_i)))
\]\[ \hspace{.02in}
= \; 
 dim_k ( e_j J(A) e_i / e_j J^2(A) e_i )
\]
This perhaps looks technical, but note that
$ e_j J(A) e_i = e_j A e_i$ when $i \neq j$ and so then is
essentially what we 
study 
in \S\ref{ss:arith1} et seq   
(and in our case the quotient
factor is even conjecturally zero, so in fact we
are already 
studying the Ext-matrix!). 
Note that the Ext-matrix defines a quiver and hence a path algebra
$kE_L(A)$.
For any finite dimensional algebra $A$, basic or otherwise, 
the Cartan decomposition matrix $C_L(A)$  is given by
\beq \label{eq:CDM}
( C_L(A) )_{ij} \; = \; dim_k Hom_A( P_j, P_i )
\eq
that is, the $i$-th row gives the number of times each simple module
occurs in $P_i$. 

\newcommand{\Mone}{\bar{{\mathcal M}}}

To pass back from the basic-algebra/homology to the full algebra we need the
dimensions of the irreducibles. For an algebra $A$ with Cartan matrix
$C_L(A)$ and a vector $v_L(A)$ giving the dimensions of the irreducible
heads of the projectives we have
\beq \label{eq:dimA}
dim(A) = v_L(A) C_L(A) v_L(A)^T
\eq

\begin{definition} \label{de:Mn1}
Let the $n \times n$ matrix $\Mone_n$  
be 
\[
\Mone_n  
  = \mat{cccccccc}
  1&\\
  1&1&\\
   &1&1&\\
   & &1&1&\\
   & & &\ddots&\ddots&\\
   & & &      &1&1&
  \tam
\]
--- we label columns left to right (and rows top to bottom)
by the ordered set $h_n$ of hook integer partitions of $n$:
$$
h_n = ((n), (n-1,1), (n-2,1^2), ..., (1^n)) .
$$
\end{definition}

We will see in \ref{th:SP}
that $\Mone_n$ is the left Cartan decomposition
matrix of $\SPe_n$
(it follows that the Ext-matrix is the same except without
the main diagonal entries).

\subsection{Linear structure of $\SPe_n$}

\medskip

A corollary of Prop.\ref{pr:SP2} is that we have an embedding of
module categories $\Gx : \func{\SPe_{n-1}}{\SPe_n}$.
In fact we can use this (together with our earlier calculations)
to determine the structure of these algebras.
Before giving the structure theorem let us recall the relevant general
theory.

\begin{lemma}{\rm (see e.g. Green \cite[\S6.2]{Green80})} \label{lem:G}
Let $A$ be an algebra and $e \in A$ an idempotent. Then
\\
(i)  the functor
$
Ae \otimes_{eAe}-:eAe-\!\!\!\!\mod \rightarrow A\!-\!\!\!\!\!\mod
$
takes a complete set of inequivalent indecomposable projective left $eAe$-modules
to a set
of inequivalent indecomposable projective $A$-modules
that is complete except for the projective covers of simple modules $L$
in which $eL=0$. (There is a corresponding right-module version.)
\\
(ii) This functor and the functor $\bar{G}_e : \func{A}{eAe}$
given by
$M \mapsto eM$ form a left-right adjoint pair.
\\
(iii) The Cartan decomposition matrix of $eAe$ embeds in that of $A$
according to the labelling of modules in (i).
\qed
\end{lemma}


\begin{figure}
  \includegraphics[width=3in]{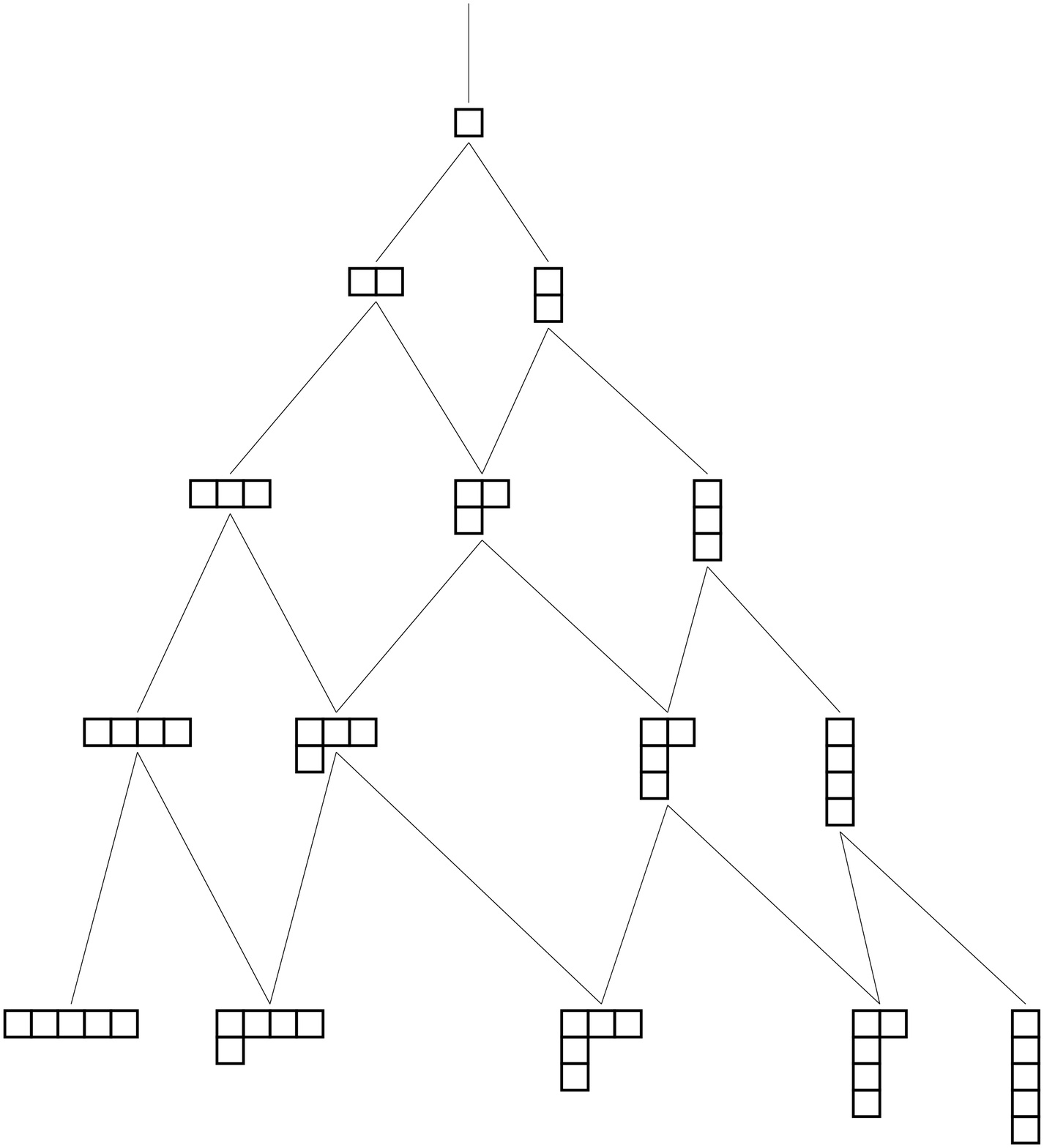}
  \caption{Young graph up to rank 5 with 22-diagrams removed. \label{fig:young1}}
\end{figure}





\begin{theorem} \label{th:SP}
  (i) Isomorphism classes of irreps 
of  $\SPe_n$   
  are naturally indexed by $h_n$.
  (Indeed $\SPe_n/rad \cong \Qq S_n/\e^1_{2,2}$ so the dimensions are
  given by the $n$-th row of Pascal's triangle
  --- see Fig.\ref{fig:young1}.) 
  \\
  (ii)  The left Cartan decomposition matrix is $\Mone_n$. 
  (Note that this determines the structure of $\SPe_n$.
  It gives the dimension as dim=$ ( 2^{n-1} + choose)/2$.)
  \\
  (iii) The image of the decomposition (\ref{eq:d1}) is complete in
  $\SPe_n$.
\end{theorem}
\proof
(i)
\ignore{{
In general given algebra $A$ and idempotent $e \in A$ then
(see e.g. \cite[\S6.2]{Green80})
the functor
$
Ae \otimes_{eAe}-:eAe-\!\!\!\!\mod \rightarrow A\!-\!\!\!\!\!\mod
$
takes a complete set of inequivalent indecomposable projective $eAe$-modules
to a set
of inequivalent indecomposable projective $A$-modules
that is complete except for the projective covers of simples $L$
in which $eL=0$. (There is a corresponding right-module version.)
}}
Consider 
Lemma~\ref{lem:G}(i).
In our case, putting $A=\SPe_n$ and $e=\chi_1$,
then by (\ref{eq:SPchi})
there is exactly one  module $L$ such that $eL=0$
(at each $n$) --- the trivial module.
Thus by Prop.\ref{pr:SP2}
$\SPe_n$ has one more class of projectives and hence irreps than
$\SPe_{n-1}$.

In particular write
$\Gx : \func{\SPe_{n-1}}{\SPe_n}$
for the functor in our case obtained using (\ref{eq:chiSP}) from
Prop.\ref{pr:SP2},
that is:
$\Gx(M) = \SPe_n \chi \otimes_{\chi \SPe_n \chi} f_\chi M$,
suppressing the
index $n$, where $f_\chi $ is as described
above. 
Then a complete set of indecomposable projectives is
\[
P^n_n = \SPe_n e^1_{(n)},
\] \[
P^n_{n-1} = \Gx(P^{n-1}_{n-1}) = \Gx(\SPe_{n-1} e^1_{(n-1)}),
\] \[
P^n_{n-2} = \Gx (\Gx( \SPe_{n-2} e^1_{(n-2)}), ...,
P^n_{n-j} = \Gx^{\circ j} (\SPe_{n-j} e^1_{(n-j)}), ...,
P^n_1 = \Gx^{\circ n-1}(k)
\]
It follows that the Cartan decomposition matrix $C(n)$ 
contains $C(n-1)$ as a submatrix, with one new row and column with the
label $n$. The new row gives the simple content of $P^n_n$.
But by (\ref{eq:sigY}) (noting Prop.\ref{th:ESSP}(ii)) this projective is simple.
Iterating, we deduce that $C(n)$ is lower-unitriangular.


Working by induction,
suppose $C(n)$ is 
of the claimed form in (ii) at level $n-1$.
Then at level $n$ we have:
\beq \label{eq:indstep1}
C(n) = \mat{c|ccccccc}
1 \\ \hline *&1 \\ * &1&1 \\ * &&1&1 \\ * &&&1&1
\\ \vdots &&&&\ddots&\ddots
\\ * &&&&       &1&1
\tam
\eq
(omitted entries 0). 
To complete the inductive step we need to compute the
$e^1_{(n)} P_{n-j}$ for each $j$.
Write $\Gx^m$ for $\Gx$ and $f^m_\chi$ for $f_\chi$ at level $m<n$, and 
note that
\[
 \Gx( \SPe_{n-1} e^1_{\lambda})
\;= \; \;
 \SPe_n \chi_1 \; \otimes_{\chi \SPe_{n}\chi} \; f_\chi ( \SPe_{n-1} e^1_{\lambda})
 \;\; =\;\;
 \SPe_n \chi_1 \; 
         \otimes_{\chi \SPe_{n}\chi}
         \;\; \chi_1 \SPe_{n-1}^{(1)} e^{(1)}_\lambda \chi_1
\] \[
 \; = \;
 \SPe_n  \chi_1 \SPe_{n-1}^{(1)} e^{(1)}_\lambda \chi_1
 \otimes_{\chi  \SPe_{n}\chi} \chi_1
\;\cong\; 
\SPe_n  \chi_1 \SPe_{n-1}^{(1)} e^{(1)}_\lambda \chi_1
\]
where
we have used that these modules are idempotently generated ideals to
apply the tensor product up to isomorphism
(and where again we use the notation from (\ref{eq:trans}),
so $\SPe_{n-1}^{(1)} $ is the 1-step translated copy of $\SPe_{n-1}$ in $\SPe_n$).
So in particular
\[
e^1_{(n)} \SPe_n \Gx( \SPe_{n-1} e^1_{(n-1)})
\;\cong\; 
e^1_{(n)} \SPe_n  \chi_1 \SPe_{n-1}^{(1)} e^{(1)}_{(n-1)} \chi_1
\;\;\subseteq \;\; e^1_{(n)} \SPe_n  \chi_1 
\]

It follows from the form of the image of $e^1_{(n)}$ in the SP
representation
(cf. \cite{Hamermesh62}, \cite[Append.B]{Martin92} and \cite{Rittenberg92})
that the dimension of $e^1_{(n)} \SPe_n  \chi_1  $  is 1, so the first $*$ is 1.
Specifically we have for example
\[
e_2 = \frac{1+p_1}{2} \mapsto \mat{c|cc|c}
1 \\ \hline
&1/2&1/2 \\
&1/2&1/2 \\ \hline
&&& 0
\tam,
\hspace{1in}
\chi \mapsto \mat{c|cc|c}
0 \\ \hline
&1&-1 \\
&0& 0 \\ \hline
&&& 1
\tam
\]
and
\[
e_3 \mapsto \frac{1}{3} \mat{c|ccc|ccc|c}
3  \\ \hline
 & 1 &1&1 \\
 & 1 &1&1 \\
 & 1 &1&1 \\ \hline
&&& & 0 &0 &0\\
&&& & 0 &0 &0\\
&&& & 0 &0 &0\\ \hline
&&&&&&&0
\tam,
\hspace{.31in}
\chi_1 \mapsto  \mat{c|ccc|ccc|c}
0  \\ \hline
 & 0 && \\
 &  &1&-1 \\
 &  &0&0 \\ \hline
&&& & 1 &-1 &0\\
&&& & 0 &0 &0\\
&&& & 0 &0 &1\\ \hline
&&&&&&&1
\tam
\]
where we have reordered the basis into fixed charge sectors,
i.e. as 111, 112, 121, 211, 122, 212, 221, 222
(the charge of a basis element is $\#(1)-\#(2)$,
where $\#(1)$ is the number of 1's
\cite{Baxter82,Martin92}).
Note from the construction that charge is conserved in SP, so each
charge sector is a submodule. We see that in each charge sector except
$(n-1,1)$ we have that either the image of $e^1_{(n)}$ is zero or the
image of $\chi_1$ is zero.
Finally in the $(n-1,1)$ sector both of these have rank 1. We deduce
that $e^1_n A \chi_1$ is 1-dimensional as required.

Similarly we have to consider
\[
 \Gx \Gx^{n-1}( \SPe_{n-2} e^1_{(n-2)})
\cong 
\SPe_n \chi_1 f_\chi f_\chi^{n-1} ( \SPe_{n-2} e^1_{(n-2)} )
\hspace{1in} \] \[
\cong 
 \SPe_n \chi_1 f_\chi  ( \chi_1 \SPe_{n-2}^{(1)} e^{(1)}_{(n-2)} \chi_1 )
= \SPe_n \chi_1 \chi_1 \chi_1^{(1)} \SPe_{n-2}^{(2)} e^{(2)}_{(n-2)}
\chi_1^{(1)} \chi_1 
\]
(NB $\chi_1^{(1)} = \chi_2$) 
giving
\[
e^1_{(n)} \SPe_n \Gx \Gx( \SPe_{n-2} e^1_{(n-2)})
\cong 
e^1_{(n)} \SPe_n f_\chi f_\chi ( \SPe_{n-2} e^1_{(n-2)} )
=
e^1_{(n)} \SPe_n \chi_1 \chi_2 ...
\]
We have, in the charge block basis, 
\[
\chi_2 \mapsto  \mat{c|ccc|ccc|c}
0  \\ \hline
 & 1 &-1& \\
 &  &0&0 \\
 &  &0&0 \\ \hline
&&& & 1 &0 &0\\
&&& & 0 &1 &-1\\
&&& & 0 &0 &0\\ \hline
&&&&&&&1
\tam,
\hspace{.31in}
\chi_1 \chi_2 \mapsto  \mat{c|ccc|ccc|c}
0  \\ \hline
 & 0 &0& \\
 &  &0&0 \\
 &  &0&0 \\ \hline
&&& & 1 &-1 &1\\
&&& & 0 &0 &0\\
&&& & 0 &0 &0\\ \hline
&&&&&&&1
\tam
\]
(in general for a nonzero entry in $\chi_1 \chi_2$ we need  basis elements with 2 in the
first and second position)
so
\[
e^1_{(n)} \SPe_n \; \chi_1 \chi_2 \; =0
\]
(remark: indeed we can verify that
$ e^1_{(n)} \; \chi_2 \chi_1 \; =0 $
holds in $\LH_n$)
so the second $*$ and indeed the other $*$s in (\ref{eq:indstep1}) are all zero.
We have verified the inductive step for (ii). 
\medskip


\noindent 
Statement (iii) 
may be deduced from (i,ii) as follows.
Note that we have $n$
isomorphism classes in the decomposition, and their multiplicities are
the dimensions of the hook irreps of $S_n$ in the natural order.
On the other hand  the $n+1$ 
charge blocks of the SP representation
are each either an irrep or contains two irreps, since each contains
one or two irreps upon restricting to $S_n$.
The first is an irrep (since dimension 1).
By the proof of (ii) the second contains the first irrep, so two
irreps in total, and the other again has the same dimension as the
corresponding $S_n$ hook representation.
Furthermore no other block contains the first irrep so this
block must be indecomposable (else the SP representation could not be
faithful, which it is by definition).
Proceeding through the blocks
then by (ii) the first $n$ of them are
a complete set of projective modules,
so each one except the first and last contains two simple
modules (`adjacent' in the hook order).
But then by the construction
of the Pascal triangle and (ii) these simple modules have the same
dimension as the corresponding $S_n$ irreps,
and (iii) follows.
\qed

\medskip

\medskip





\section{On representation theory of \texorpdfstring{$LH_n$}{LHn}} \label{ss:lowrank}


\newcommand{\chim}[1]{\chi^{(#1)}} 
\newcommand{\MMMMt}{ \mat{cccccccccc} 1&1 \\ &1&1 \\ &&1&1 \\ &&&1&1 \\ &&&&\ddots \tam  }
\newcommand{\MMMM}{ \mat{cccccccccc} 1& \\ 1&1 \\ &1&1 \\ &&1&1 \\ &&&\ddots&\ddots \tam  }
\newcommand{\MMMMM}{ \mat{cccccccccc}
1&1&1&1&1&... \\ &1&1&1&1&... \\ &&1&1&1&... \\ &&&1&1&... \\ &&&&\ddots \tam  }

Combining \eqref{eq:Pascn} with \eqref{eq:dimA} and Th.\ref{th:SP} we have 
\[
dim(\SPe_n(t\neq 1)) = v_1 \UU^{n-1} \MMMM \left( v_1 \UU^{n-1}  \right)^T 
=\binom{2n-1}{n-1}=\frac{1}{2}\binom{2n}{n}
\]
--- NB we have used the obvious `global' limit of all the Cartan matrices
(it is a coincidence that this and the $\UU$ matrix are similar).

\newcommand{\Diag}{Diag}

Given a vector $v$ we write $\Diag(v)$ for the diagonal matrix with
$v$ down the diagonal.
Let $p^n$ be the vector with the $n$-th row of Pascal's triangle as
the entries, thus for example $p^4 = (1,3,3,1)$. 
We have
\[
\MMM^p_n \; :=\; \Diag(p^n)\; \bar{\MMM}_n \; \Diag(p^n)
\]
(examples are given in \eqref{eq:Mpn} below)
and the dimension is the sum of all the entries. The closed form
follows readily from this.
Also from Th.\ref{th:SP} we have:


\begin{corollary} 
  \label{con:morita}
For $t \neq 1$ 
  the Morita class of $\SPe_n$ is of the path algebra with $A_n$ quiver
(directed $1 \rightarrow 2 \rightarrow \cdots \rightarrow n$)
  and relations given by vanishing of all proper paths of length 2.
  In particular the radical-squared vanishes.
\ignore{{
  \\  The structure of the quotient
  $\SPe_n/\chim{j+1}$  
  is given by the $j\times j$ truncation of
\ppm{ fix: } $\MMM_n$.
\todopm{(NB havn't checked line up of indices yet.)} 
}}
\end{corollary} 

\ignore{{
\medskip

{\em Remark}. Anyway it seems that something rather nice is going on
here in the crossover between the low-d topology of the loop braid and
the algebra!... More on this later.

\bigskip
\medskip

Observe that $\rrb{GB}$ contains, 
by restriction, the
$(n-1,1)$-Young representation of the corresponding symmetric group.
Thus, restricting to the corresponding subalgebra, $\rrb{GB}$ decomposes as
$n=(n-1)+1$ --- in this case a direct sum decomposition.
Noting that $\rrb{GB}$ always has a submodule corresponding to the
trivial representation, we see that it also decomposes as $n=(n-1)+1$
over the full algebra --- however we do not yet know if this is a
direct sum decomposition. 
For example for $n=3$ the 3-dimensional $\rrb{GB}$ decomposes as
$
3=1+2
$
on the subalgebra corresponding to the $\rr{i}$s.
This bounds below the decomposition of the 3 in the representation of the full
algebra.
However, 
(for example by Magma)  
the irreducible dimensions of $LH_3$ are
$1+2+1$ (obvious shorthand).
We deduce that $\rrb{GB}$ `contains' the
2-dimensional irrep as a subquotient.
It is also evident that $\rrb{GB}$ contains the trivial irrep as a
submodule, as already noted.
But, as we will see below,
$\rrb{GB}$  is not  semisimple, so these components do not form a
direct sum.

As an aside, for $n=2$ we have
\[
\rrb{GB} (\sig{1}\rr{1}-\rr{1}\sig{1})
=
\mat{ccc}
t-1 & 1-t \\
t-1 & 1-t
\tam
= (t-1) \mat{ccc} 1&-1 \\ 1 & -1 \tam
\]
so the representation here is not a direct sum
(this is an element of the radical - see \ref{} later - with non-zero image). 

 It follows from the above that we can compute the characters of the
irreducible representations of $\LH_n$ in case $n=3$.
Hence we can determine the irreducible content of the regular
representation. 
}}

\medskip






\subsection{Properties determined from Th.\ref{th:SP} and direct calculation in low rank}
\label{ss:low}
%



Our results for $\LH_n$ may be neatly given as follows.
Firstly,
\begin{proposition}
For $t^2 \neq 1$ and $n<8$, $\LH_n \cong SP_n$.
\end{proposition}
\proof Here we can compute dimensions directly, which saturates the bound
on the kernel. \qed
\begin{conjecture} \label{con1x}
For $t^2 \neq 1$, $\LH_n \cong SP_n$.
\end{conjecture}
\ignore{{
Then we note that {\em the conjecture is confirmed up to $n=7$}
by direct calculation of dimension.
\\
Unpacking this gives a nice way to read off properties from the SP
structure Theorem.
}}



A summary of what we learn for the algebra dimensions,
and  irreducible reps, of $\LH_n$ is given by the following tables.

\beq \label{eq:Pascal}
\begin{array}{c|c|c|c|c|ccccccccccccccc}
&t=1   & t=-1   &t^2\neq1&t\neq 1 & irreps: \hspace{-.3in}&&&\;labels \hspace{-.32in} \\
   &    &     & &  & /dims  \hspace{-.2in} \\
  n&dim& dim&dim&ssdim&  -6 &-5&-4&-3&-2&-1&0&1&2&3&4&5 &6 \\
  \hline
1&1   &1&   1&  1&&&  &&&&1 \\
2&3   &3&   3&  2&&&  &&&1&&1 \\
3&15  &11&  10&  6&&&  &&1&&2&&1 \\
4&114 &42&  35& 20&&&  &1&&3&&3&&1 \\
5&1170 &163 &126& 70&&&  1&&4&&6&&4&&1 \\
6&15570 &638 &462&252& &1&&5&&10&&10&&5&&1 \\
7&   & 2510&1716  & 924  & 1&&6&&15 & &20&&15&&6& & 1 
\end{array}
\eq
The irrep labels here are given by $(n-i,1^i) \mapsto n-2i-1$.

\newcommand{\MMMMMT}{ \mat{cccccccccccc}
   1&&
\\ 1&1&
\\ 1&1&1&
\\ 1&1&1&1&
\\ \vdots&\vdots&\vdots&\vdots&\ddots
\tam  }

Combining \eqref{eq:Pascn} with \eqref{eq:dimA}, Th.\ref{th:SP} and \ref{con1x} we
have the conjecture
\[
dim(\LH_n(t^2\neq 1)) = v_1 \UU^{n-1} \MMMM \left( v_1 \UU^{n-1}  \right)^T 
=\binom{2n-1}{n-1}=\frac{1}{2}\binom{2n}{n}
\]
For $t=-1$ we note that $\SPe_n$ is generally a proper quotient of
$\LH_n$,
and that $\LH_n$ has larger radical (the square does not vanish).
We define the semiinfinite matrix:
\[
C(\LH_{}(t=- 1)) =  \MMMMMT
\]
and conjecture that the  Cartan matrix  $C(\LH_{n}(t=- 1))$
is this truncated at $n \times n$
(i.e. the quiver is the same as the generic case, but without quotient
relations);
and thus we conjecture
\[
dim(\LH_n(t=- 1)) \; =\; v_1 \UU^{n-1} \MMMMMT \left( v_1 \UU^{n-1}  \right)^T 
 = \frac{n^2 + \binom{2n-2}{n-1}}{2}
\]
(cf. OEIS A032443). Note that our calculations verify this for $n\leq 7$.

For $t=1$ we see that $LH_n(t=1)$ has semisimple quotient at least as
big as $\mathbb{C}S_n$, which is of dimension $n!$.
Indeed, in this case the quotient by the relation $\sigma_i=\rho_i$ is
precisely $\mathbb{C}S_n$, since in this case $\sigma_i^2=1$.  For
$n\leq 4$ we have computationally verified that the semisimple
subalgebra of $LH_n(t=1)$ is precisely $\mathbb{C}S_n$, and we
conjecture  
that this is the case for all $n$.  The Jacobson radical grows quite
quickly however,
and we do not have a conjecture on the general structure.

\ignore{{
Define a semiinfinite matrix by
\[
\MMM = \mat{cccccccccc} 1&1 \\ &1&1 \\ &&1&1 \\ &&&1&1 \\ &&&&\ddots \tam
\]
}}

\ignore{{
Noting the apparent Pascal triangle here, let us set some labels.
We label the positions in the triangle according to the `charge', the 
number of steps more to the right than the left.
Thus a step sequence
$+----$ gives $-(n-2)=-3$ and so on.
Alternatively we label by `weight' $\lambda=(R,L)$ where $R$ (resp. $L$) is the number of
right (resp. left) steps,
and our example gives $(1,n-1)$.
The connection from weight to charge is $(R,L) \mapsto R-L$. 

We have in mind that $-n $ (reached by the step sequence $ ---\cdots -$)
labels the `other' 1d representation and
\[
+n \leadsto +++\cdots +
\]
(step right, step right, ...) labels the `trivial' representation.
\todopm{-or other way round.}
Intrinsic characterisations of the other irreducibles have yet to be
given. 
}}

Observe that the numbers in \eqref{eq:Pascal} follow the conjectured patterns.
Since the vector $v_1$ has finite support the nominally infinite sums
above are all finite.
To inspect the supported part, in the generic case
consider matrices $\MMM^p_n$ ($n=2,3,4,5$):
\ignore{{
\[
\mat{ccc} 1&1 \\ & 1 \tam, \qquad
\mat{cccc} 1&2 \\ & 2^2 & 2 \\ && 1 \tam ,
\qquad
\mat{ccccc} 1&3 \\ & 3^3 & 9 \\ && 3^3 & 3 \\ &&& 1 \tam ,
\qquad
\mat{ccccccc} 1&4 \\ & 4^2 &24 \\ && 6^2 & 24 \\ &&& 4^2 & 4 \\ &&&&1  \tam
\]
}}
\beq \label{eq:Mpn}
\mat{ccc} 1& \\ 1 & 1 \tam, \qquad
\mat{cccc} 1& \\ 2& 2^2 &  \\ &2& 1 \tam ,
\qquad
\mat{ccccc} 1& \\ 3 & 3^3 &  \\ &9& 3^3 &  \\ &&3& 1 \tam ,
\qquad
\mat{ccccccc} 1& \\ 4& 4^2 & \\ &24& 6^2 &  \\ &&24& 4^2 &  \\ &&&4&1  \tam
\eq
Here the semisimple dimension is given by the sum down the diagonal
and the radical dimension is given by the sum in the 
off-diagonal.
\ignore{{
The next 
is $n=6$:
\beq
\MMM^p_6 = 
\mat{cccccccc}
1&5 \\ &5^2&50 \\ &&10^2&100 \\ &&&10^2&50\\ &&&&5^2&5 \\&&&&&1
\tam^T
\qquad \mbox{ giving dim=}252+210=462
\eq
}}
\ignore{{  
Of course the sum down the diagonal is also a central binomial
coefficient (thus $(2m)!/(m!)^2$ where $m=n-1$).
The other sum  also has a nice closed form expression.
}}

For $t=-1$:
\[
\mat{ccc} 1& \\ 1 & 1 \tam, \qquad
\mat{cccc} 1& \\ 2& 2^2 &  \\ 1&2& 1 \tam ,
\qquad
\mat{ccccc} 1& \\ 3 & 3^3 &  \\ 3 &9& 3^3 &  \\1 &3&3& 1 \tam ,
\qquad
\mat{ccccccc} 1& \\ 4& 4^2 & \\ 6&24& 6^2 &  \\ 4&16&24& 4^2 &  \\ 1&4&6&4&1  \tam
\]


\subsection{On $\chi$ elements}

Let us define
\beq
\chim{m+1} = (\s_1 -\p_1)(\s_2 -\p_2)...(\s_m -\p_m) ,
\eq
understood
as an element in $LH_n$ with $n>m$.
Thus in particular $\chim{2} = \chi_1$. 
Similarly
for sequence $X=(x_1, x_2, ..., x_k)$ define
\beq
\chim{X}  = (\s_{x_1} -\p_{x_1})(\s_{x_2} -\p_{x_2})...(\s_{x_k} -\p_{x_k}) ,
\eq
and
\beq
\chim{m+1}_-  = (\s_m -\p_m)(\s_{m-1} -\p_{m-1})...(\s_2 -\p_2)(\s_1 -\p_1) ,
\eq
It is easy to verify that if $X$ is non-increasing then
$\chim{X} \chim{X}  =  (1-t)^k \chim{X}$.
Thus (for $t \neq 1$) the non-increasing cases can all be normalised
as idempotents.
However it is also easy to check that no increasing case can.
(A nice illustration of the `chirality' present in the defining relations.)

Observe that imposing the relation $\s_1 = \p_1$ in $LH_n$ forces $\s_1 =1$,
unless $t=1$.
Thus the quotient algebra
\beq \label{eq:q2triv}
LH_n / \chim{2} \; \cong k \hspace{1in} t \neq 1
\eq
i.e. only the trivial, or label $\lambda=+n$, irrep survives.
And the same holds for $\SPe_n$.
The following has been checked up to rank 5.

\ignore{{

\ppm{(This para will be deleted once assimilated elsewhere.)}
Observe also that truncating
the matrix $\MMM^p_n$ 
in the top left $j \times j$ corner we
get, varying $j$ appropriately,
the dimensions of all the quotients we have seen so far.
For example $j=2$ gives dim $A_n$
--- this is proved below
(NB it is conjectured below that this is also the quotient by
$(s_1 -p_1)(s_2 -p_2)$).
But also in cases checked so far 
$j=3$ gives the dim of the quotient by
$(s_1 -p_1)(s_2 -p_2)(s_3 -p_3)$,
for example with $n=4$ we get 31 and with $n=5$ we get 81
(all as per Eric's magma).

}}

\medskip

\begin{conjecture}
  \label{con:morita2}
  The structure of the quotient
  $LH_n/\chim{j+1}$
  is given by the $j\times j$ truncation of
$\MMM^p_n$.
\end{conjecture}






\section{Discussion and avenues for future work}

Above we give answers to the main structural questions for $\SPe_n$ and $\LH_n$.
But exploration of generalisations is also well-motivated,
since these algebras
(even taken together with the constructions discussed in \cite{Kadar-Martin-Rowell-Wang:2016})
cover a relatively small quotient inside
$Rep(LB_n)$.
With this in mind, there
are a number of other questions worth addressing around $\SPe_n$ and
$\LH_n$,
offering clues on generalisation, and hence 
towards
understanding
more of the structure of the group algebra.  
Remark \ref{remark:other quotients} suggests that for most values of $t$ we obtain larger finite dimensional quotients
by eliminating one of the local relations (\ref{eq:R1i}) or (\ref{eq:R1ii}).  Computational experiments suggest that for $t=0$ eliminating (\ref{eq:R1ii}) yields infinite dimensional algebras.  This parameter-dependence should be further explored.

In light of the results of \cite{Reutter:preprint} the non-semisimplicity of $LH_n$ is an important feature, rather than a shortcoming.  Extracting topological information from the non-semisimple part requires some further work, as Markov traces typically `see' the semisimple part.  Another aspect of our work is the (conjectural) localisation of the regular representation of $LH_n$.  It is worth pointing out that localisations of \emph{unitary} sequences of $B_n$ representations are relatively rare, conjecturally corresponding to representations with finite braid group image \cite{Rowell-Wang:2012,GHRIMRN}.  Since $LH_n$ is non-semisimple and hence non-unitary this does not contradict this conjectural relationship, but gives us some hope that localisations are possible for other parameter choices and other quotients.




The quotient of $LB_n$ by the relation ${\sigma_i}^2=1$ is a
potentially interesting infinite group, which we call the mixed double
symmetric group $MDS_n$.  The reason for this nomenclature is that
$MDS_n$ is a quotient of the free product of two copies of the
symmetric group.  In particular, $MDS_n$ surjects onto
$S_n$ by ${\sigma_i}\rightarrow \rho_i$.  It is of special interest
here as $LH_n(1)$ is a quotient of $\mathbb{Z}[MDS_n]$.
We expect it could be of quite general interest.

In \cite{Kadar-Martin-Rowell-Wang:2016} constructions are developed
based on BMW algebras, but still starting from `classical'
precepts. It would be very interesting to meld the super-Burau-Rittenberg
construction to the KMRW construction.  For example, one might try to use cubic local (eigenvalue) relations among the generators $\rho_i,\sigma_i$ to obtain finite dimensional quotients, possibly inspired by the relations satisfied by a subsequence of $LB_n$ lifts of BMW algebra representations.

\ignore{{
something about higher Rittenbergs / tensor reps for $t=-1$ / blob
/
Chimera: our localisation functor places attention on $\chi_1 = \chim{2}$.
/
see old appendix
/
\ppm{actually put something here!?... or not....  :-)} 
}}




\vspace{1cm} \appendix

\section{Preparatory arithmetic and notation for left ideals}



\ignore{{
\begin{figure}
  \includegraphics[width=3in]{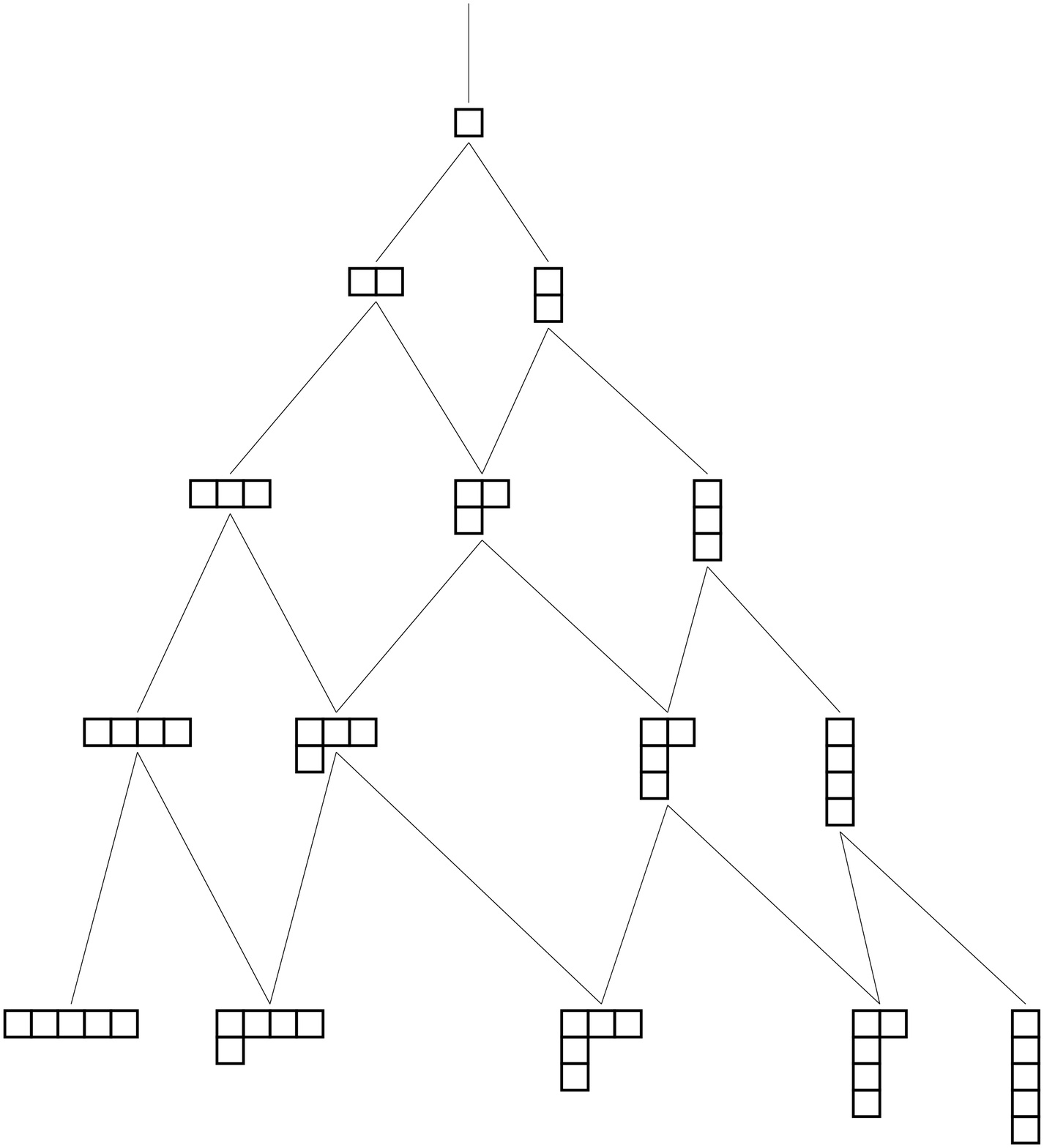}
  \caption{Young graph up to rank 5 with 22-diagrams removed. \label{fig:young1}}
\end{figure}

\bigskip

}}


\subsection{
  Symmetric group and Hecke algebra arithmetic}\label{ss:SnIds}

Recall Young's (anti-)symmetrizers in $kS_n$.
Unnormalised in $\Zz S_n$ they are:
\beq
Y_{\pm}^n  = \sum_{g \in S_n} (\pm 1)^{len(g)} g
\eq
where $len(g)$ is the usual Coxeter length function.
If $k$ has characteristic 0 then $kS_n$ is semisimple and these
elements are simply the (unnormalised) idempotents corresponding to
the trivial and alternating representations respectively.
Note that exactly the same classical construction works for the Hecke algebra
over any field where it is semisimple.
(The corresponding idempotents are sometimes called Jones--Wenzl
projectors.)
Specifically (see e.g. \cite[\S9B]{CurtisReiner81})
\[
X_{\pm}^n  = \sum_{g \in S_n} (-\lambda_\mp )^{-len(g)} T_g
\hspace{1in}
\mbox{i.e. }
X_-^2 = 1-\sig{1} , \;\; X_+^2 = 1+t^{-1} \sig{1} , ...
\]
where for us $\lambda_- =-t$ and $\lambda_+ = 1$
(the apparent flip of labels is just because we use non-Lusztig scaling),
and $T_g$ is the product of
generators obtained by writing $g$ in reduced form then applying
$\rr\ii \mapsto \sigT\ii$.

Working in $kS_{n+m}$  
we understand $Y^n_{\pm}$
and translates such as $Y^{n(1)}_{+}$ in the obvious way.
Note then that we have many identities like
\beq \label{eq:Yid1}
Y^2_{+} Y^n_{+} = 2 Y^n_{+},
  \hspace{1in}   Y^{a(1)} Y^n_{+} = a! Y^n_{+} \quad (a<n)
\eq

Recall  $\Lambda_n $ denotes the set of integer partitions of $n$.
Over the rational field we have a decomposition of $1 \in kS_n$ into
primitive central idempotents
\beq
1 = \sum_{\lambda\in\Lambda_n} \epsilon_\lambda
\eq
where  
each $\epsilon_\lambda$ is a known unique element
(see e.g. Cohn \cite[\S7.6]{Cohn82II} or Curtis--Reiner \cite{CurtisReiner81} for
gentle expositions).
There is a further (not generally unique)
decomposition of each $\epsilon_\lambda$ into
primitive orthogonal idempotents
\beq
\epsilon_\lambda = \sum_{i=1}^{dim_\lambda} \e_\lambda^i
\eq
where $dim_\lambda$ is the number of walks from the root to $\lambda$
on the directed Young graph.
The elements $\e_\lambda^i$ are conjugate to each other.
The elements $\e_\lambda^i$ are
not uniquely defined in general. 
Two possible constructions
of one for each $\lambda$ are
exemplified pictorially by (case $\lambda=442$)
\beq \label{eq:choicee}
  {e}_\lambda^{1} = \;\;\;
  c_\lambda \;\raisebox{-.321in}{\includegraphics[width=2.095cm]{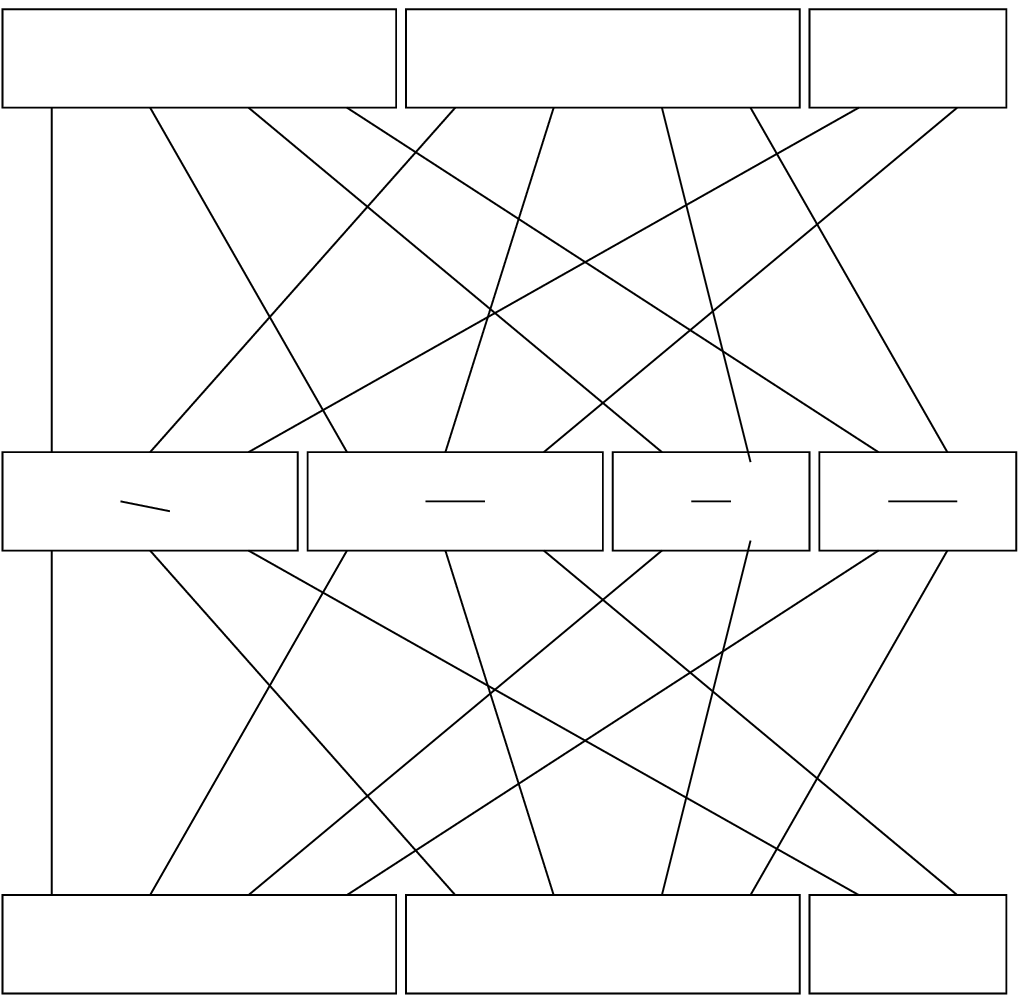}}
, \hspace{1in}
  \hat{e}_\lambda^{1} = \;\;\;
  c_\lambda \;\raisebox{-.241in}{\includegraphics[width=2.095cm]{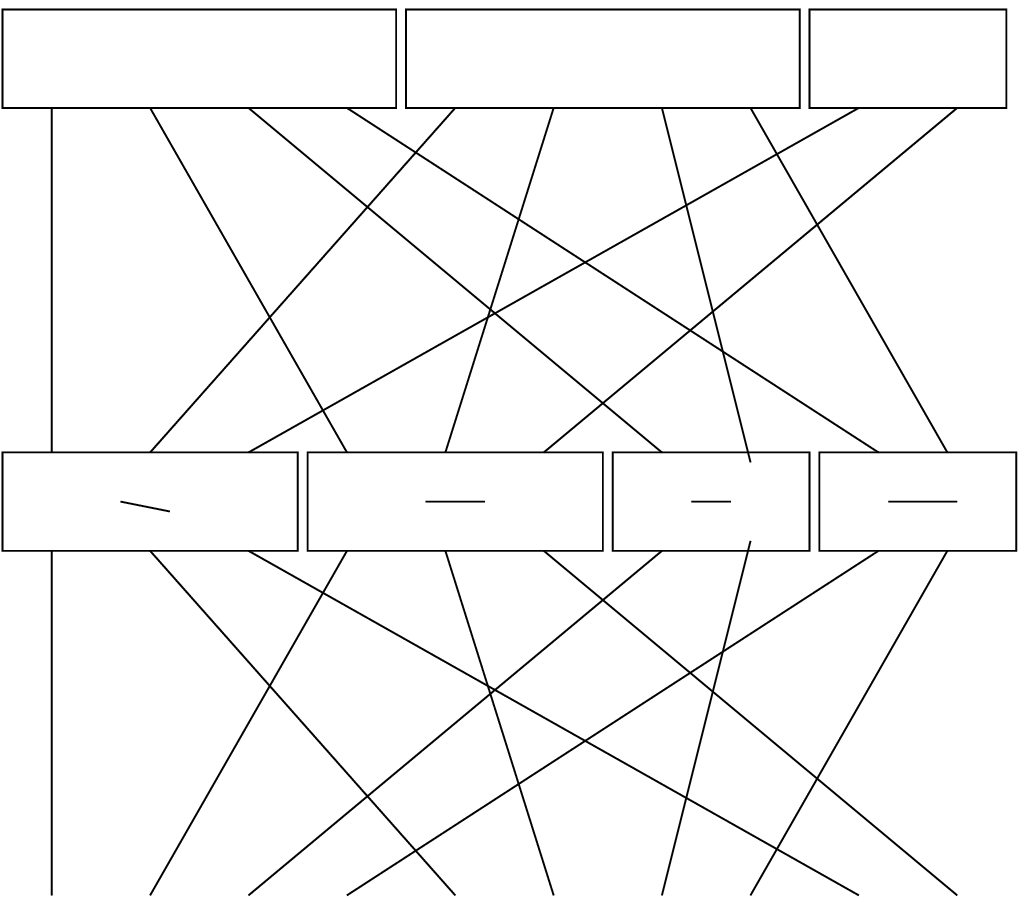}}
\eq
where an undecorated box is a symmetrizer and a $-$-decorated box an
antisymmetrizer,
and the factor $c_\lambda$ is just a scalar.
(NB For the moment we write $ e_\lambda^{1}$ instead of $ \e_\lambda^{1}$
for this specific choice.)
In particular though,
$\e_{(n)}^1$ {\em is} unique:
$ \e_{(n)}^1 = \frac{1}{n!} Y_+^n$.
(The whole story lifts to the Hecke case --- see e.g. \cite{Martin:M1}
for a gentle exposition.)

\medskip

An idempotent decomposition of 1 in a subalgebra $B$ of an algebra $A$
is of course a decomposition in $A$.
Thus in particular we can take an idempotent in $k S_n$ and consider
it as an idempotent in $k S_{n+1}$ by the inclusion that is natural
from the presentation ($p_i \mapsto p_i$).
Understanding $\e_{\lambda}^1 $ with $\lambda\vdash n$ in $ k S_{n+1}$ in this way, 
a useful property in our $k=\Cc$ case will be
\beq \label{eq:idres}
\e_\lambda^j = \sum_{\mu \in  \lambda+} \e_{\mu}^{'}
\eq
where $\lambda+$ denotes the set of partitions obtained from $\lambda$
by adding a box, and the prime indicates that we identify this
idempotent only up to equivalence.
(Various proofs exist. For example note that the existence of such a decomposition
follows from the induction rules for $S_n \hookrightarrow S_{n+1}$.)
For example
\[
\e_{(2,2)}^1 = \e_{(3,2)}^{'} + \e_{(2,2,1)}^{'}
\]


\newcommand{\germ}{\mathfrak}  
\newcommand{\cprime}{'} 

\bibliography{bib/Loop_Hecke.bib,bib/add-Deguchi.bib}{}
\bibliographystyle{alpha}

\noindent
CD, PM: School of Mathematics, University of Leeds, Leeds, UK.
\\ ECR: Department of Mathematics, Texas A \& M University, Texas, USA.
\end{document}